\documentclass{article}

\newtheorem{theorem}{Theorem}
\newtheorem{lemma}[theorem]{Lemma}
\newtheorem{definition}{Definition}  
\newtheorem{proposition}[theorem]{Proposition}

\newcommand\M{{\mathcal M}}
\newcommand\N{{\mathcal N}}
\newcommand\J{{\N}}
\newcommand\A{{\mathbf A}}
\newcommand\D{{\mathcal D}}
\newcommand\Tr{{Tr_{\M,c}}}
\newcommand\Trro{{Tr_{\M,c,\ro}}}
\newcommand\Trk{{Tr_{\M,c}^k}}
\newcommand\T{{ ~^1Tr_{\M,c}}}

\newcommand\ro{\rho}

\newcommand\Smallskip{

\smallskip\noindent}

\newcommand\Prec{\prec_\tau}
\newcommand\Preceq{\preceq_\tau}
\newcommand\Precs{\prec^s_\tau}
\newcommand\K{KB_\tau}

\newcommand\m{Let $\M$ be a structure}

\newcommand\Benum{\begin{enumerate}}
\newcommand\Eenum{\end{enumerate}\vspace{-3pt}}
\newcommand\Item{\item\vspace{-3pt}}

\newcommand\Theorem[1]{\begin{theorem}\label{#1}
}
\newcommand\EndTheorem{\end{theorem}}

\newcommand\Proposition[1]{\begin{proposition}\label{#1} }
\newcommand\EndProposition{\end{proposition}}

\newcommand\Lemma[1]{\begin{lemma}\label{#1}
}
\newcommand\EndLemma{\end{lemma}}

\newcommand\Definition[1]{\begin{definition}\label{#1}}
\newcommand\Enddefinition{\end{definition}}

\newcommand\Corollary[1]{\noindent{\bf Corollary #1}\it}
\newcommand\Endcorollary{\rm}

\newcommand\Proof{\noindent {\tt Proof. }}
\newcommand\Endproof{{\tt QED}}

\usepackage{fullpage}
\begin{document}

\title{On the hierarchies of universal predicates}
\author{Pavel Hrube\v{s}\footnote{{\it Mathematical Institute, Academy of Sciences of the Czech
  Republic, Prague, Czech republic  }}}
\date{}
\maketitle 

\vspace{-18pt}
\begin{abstract}
We investigate a hierarchy of arithmetical structures
obtained by transfinite addition of a canonic universal
predicate, where the canonic universal predicate for $\M$ is defined as a
minimum universal predicate for $\M$ in terms of definability. We
determine the upper bound of the hierarchy and give a characterisation
for the sets definable in the hierarchy.
\end{abstract}


\section{Introduction}\label{introduction}
In his fundamental works, Bernard Bolzano develops the idea that all
natural languages are approximations of a single \emph{universal
  language;} a language in which we can describe anything that exists or
that could exist. The idea was not new at Bolzano's time and it still
persists, at least as a call in the hearts of
logicians. However, if there is a general conclusion which can be drawn
from the results related to the G\"odel theorem then it is this: 
\emph{there is no universal language.} 
For any language $L$, there exists at least one thing that cannot be
fully described in the language: the semantics of $L$ itself.  
 In this
paper, we shall be working with first-order languages interpreted over
natural numbers, but we believe that some of the results are valid in
general. Let $\M$ be a first-order structure with a countable language
interpreted over $\omega$. According to the well-known theorem of
Tarski, the truth predicate for $\M$ (i.e., the set of G\"odel
numbers of sentences true in $\M$) is not definable in $\M$. 
It must
be noted that this proposition does not merely assert that there
is a set which cannot be defined in $\M$, but it gives an example of
such a set; and moreover, the set, as a description of the semantics
of $\M$, is presupposed in the structure
$\M$ itself.  If we take a structure $\M_0$ and the truth
predicate $T_0$ for $\M_0$ then the structure $\M_1:=M_0+T_0$ will be
stronger then $\M_0$. Moreover, it is an extension in a sense presupposed already
in the structure $\M_0$. Similarly we can define
$\M_2:=\M_1+T_1$ etc. and we can even imagine that we iterate the
process transfinitely and obtain an infinite hierarchy of 
structures $\{M_\alpha\}$. The structures in the hierarchy are natural extensions of $\M_0$
and it makes sense to ask what are the properties of such a
hierarchy, which sets are definable at some stage of the sequence
etc. 
The notion of the hierarchy, however, is a reminiscence of the idea of the
universal language, and it must inevitably lead into difficulties. 
The first and main problem is that the notion 'truth for $\M$' is not
determined uniquely. More generally, we want the structure
$\M_{\alpha+1}$ to be obtained as \emph{'$\M_{\alpha}$+'the description
  of the semantics of $\M_\alpha$'}.  However, the notion \emph{'the
description of the semantics of $\M_\alpha$' } is not
unambivalent, as there may exist infinitely many sets which may be
said to describe the semantics of $\M_\alpha$. It is an obvious move
to try to chose a
particular, \emph{canonic,} description of semantics of $\M_\alpha$
and define the hierarchy in terms of adding the canonic
description. Two alternative definitions of such a description
will be given below under the headings \emph{canonic universal
  predicate} and \emph{proper canonic universal predicate}. Of course,
we must
then answer the question whether such a canonic choice is
possible, i.e., we must determine whether a (proper) canonic universal
predicate for a given $\M$ exists, and this problem will form the major part of
the present paper.

\bigskip
For rather technical reasons (explained on page \pageref{defIc})  the truth predicate itself
is not exactly suitable for the purpose of defining a hierarchy, and
we shall thus define the hierarchy in a related but different
way. Moreover,
 we shall investigate two kinds of hierarchies, one
obtained using \it a proper universal predicate, \rm and the other
using \it a universal predicate for $\M$
\rm. Let ${\mathcal D(\M)}\subseteq {\mathcal P}(\omega)$ denote the
set of all (first-order) definable sets in $\M$,
 then\footnote{If $X\subseteq \omega^2$ and $n\in \omega$, we set
$X(n,.):=\{m\in\omega;X(n,m)\}$.} \Benum
    \Item $G\subseteq \omega^2$ is \emph{ a universal predicate
 for $\M$ } iff ${\mathcal D}({\cal M})\subseteq \{G(n,.);n\in \omega\}$
\Item $P\subseteq \omega^2$ is \emph{ a proper universal predicate for $\M$
 }iff ${\mathcal D}({\cal M})=\{P(n,.);n\in \omega\}$.
\Eenum
\noindent
 Evidently, a (proper) universal
predicate for $\M$ is not definable in $\M$; but it must be observed that neither $P$ nor $G$ are unique for a given
$\M$. Moreover, $P$ or $G$ can be chosen in such a way that the
structures $\M+P$ and $\M+G$ can have an arbitrary strength. In order
to avoid the problem we introduce a \it canonic (proper) universal
predicate \rm as a minimum (proper) universal predicate in the
following sense

\smallskip\label{canonic}
$P_0$ is \emph{ a canonic (proper) universal predicate for $\M$} iff
\Benum
\Item  $P_0$ is a (proper) universal predicate for $\M$ and 
\Item for every (proper) universal predicate $P$ for $\M$, $P_0$ is definable in $\M+P$.
		\Eenum

The obvious question is whether the
canonic universal or canonic proper universal predicates exist. The answer,
which is partially given in this paper, is non-trivial: there
are structures which have a canonic (proper) universal predicate and
there are countable structures which do not. Further, the two concepts are not equivalent and there are structures which possess a canonic proper universal predicate but do not have a canonic universal predicate. 

For a given countable structure $\cal N$ we shall define the \emph{Tarski hierarchy{\footnote{We use the name of Tarski
  because he was the first one to state the undefinability of truth. I
am not aware that he would ever attempt to iterate the process of
adding the truth predicate finitely or transfinitely.}} over $\cal N$ }
to be a sequence of structures $\{M_\alpha\}_{\alpha\leq\lambda{(\cal N})}$ such that \Benum
   \Item ${\mathcal M}_0=\mathcal N$, 
   \Item for every $\alpha<\lambda({\cal N})$, $\cal M_\alpha$ has a canonic universal predicate and ${\mathcal M}_{\alpha+1}={\mathcal M}_\alpha\cup \{P_\alpha\}$,
   \Item $\mathcal M_\beta=\bigcup_{\beta<\alpha} M_\alpha$ for every
         limit ordinal 
$\beta\leq\lambda({\cal N})$, 
   \Item $\lambda(\N)$ is the maximum ordinal satisfying 1)-3).
\Eenum
An analogous hierarchy obtained by
replacing the notion {\it canonic universal predicate } by that of
{\it canonic proper universal predicate }  and $\lambda(\N)$ by
$\lambda^p(\N)$ will be called \emph{
 proper Tarski hierarchy over $\N$}.\label{Tarski}
In essence, proper Tarski hierarchy can be viewed as a sequence of
truth predicates. 

The key characteristic of the Tarski hierarchy is the ordinal $\lambda(\cal
\N)$. A priori, we know that $\lambda({\cal N})\leq\aleph_1$, as an
uncountable structure cannot have a universal predicate. If
$\lambda({\cal N})=\aleph_1$ then every countable structure obtained by
the process of adding a canonic universal predicate does possess a
canonic proper universal predicate. If, on the other hand, we have
$\lambda({\cal N})<\aleph_1$ then the structure $\M_\lambda$ does not
have a canonic universal predicate. 

 The basic properties of the
proper Tarski hierarchy can be obtained from [1]. The authors define the hierarchy as a sequence of Turing
degrees $\{H_\alpha\}_{\alpha<\xi}$. On the isolated steps they take for
$H_{\alpha+1}$ simply the Turing jump for $H_{\alpha}$ and the
minimum-proper-universal-predicate question enters on the limit steps. But in
essence, their definition is equivalent to the notion of Tarski
hierarchy adopted here.\footnote{ Let $\{H_\alpha\}$ be the sequence of
  Turing degrees in the sense of [1] and let $\{\M_\alpha\}$ be a
  proper Tarski hierarchy over $\A$. The interrelation between the
  hierarchies is based on the following fact: if
  $\M_\beta\dot=H_\gamma$, where $\dot=$ is understood in terms of
  Turing reducibility, then $\M_{\beta+1}\dot=H_{\gamma+\omega}$. Since every element in $H_\alpha$ is
  finite, there is no counterpart in $\{H_\alpha\}$ to $\M_\beta$
  if $\beta$ is a limit. But it is trivial to find some $\gamma$ such
  that $\M_{\beta+1}\dot=H_{\gamma}$.  } In particular, using the techniques developed
in their paper, it can be shown that for every countable $\cal N$ the
ordinal $\lambda^p(\cal N)$ is a countable limit ordinal.  An alternative
approach to proper Tarski hierarchy and a
comparison of Tarski and proper Tarski hierarchy can be found in [3].  

In this paper  
we shall be concerned mainly with the Tarski hierarchy.  
We shall prove the following main results

\Theorem{Tha}
Let $\{\M_\alpha\}_{\alpha\leq\lambda(\cal N)}$ be a 
Tarski hierarchy over a $L$-finite structure $\cal N$. Then $\lambda(\cal
N)$ is countable. Furthermore, $\lambda({\cal N})=Ord(\N)$, the first
undefinable (i.e., non-recursive) ordinal in $\cal N$, and the structure ${\cal M}_{\lambda(\N)}$ is the
minimal structure containing all sets implicitly definable in $\N$.
\EndTheorem

 We shall note that the structure ${\cal M}_{\lambda(\N)}$ is at the
 same time the smallest structure containing all sets
 $\Delta_1^1$-definable in $\cal N$, ie. it coincides with the sets
 hyperarithmetical in $\N$. 

\Theorem{Thb}
Let $\{\M_\alpha\}_{\alpha\leq\lambda(\cal N)}$ be a Tarski hierarchy 
over a $L$-finite structure $\cal N$ and let 
$\{\M^p_\alpha\}_{\alpha\leq\lambda^p(\cal N)}$ be a proper Tarski hierarchy over
$\cal N$. 
Then $\lambda^p(\cal N)>\lambda(\cal N)$ and for every
$\alpha\leq\lambda(\cal N)$,  $M_\alpha^p\sim M_\alpha$. Hence
the structure $M_{\lambda(\cal N)}$ does not have a canonic universal predicate but
  does have a canonic proper universal predicate.
\EndTheorem

The part of Theorem~\ref{Tha} asserting that  $\lambda(\N)\geq
Ord(\N)$ is proved as Theorem
\ref{LemmaIIIoq}. That the structure $\M_{Ord(\N)}$ is the smallest
structure containing all implicitly definable sets in $\N$ is claimed in
Theorem~\ref{Thx} and proved on page \pageref{proofx}. Finally, the
fact that $\M_{Ord(\N)}$ does
not have a canonic universal predicate and hence $\lambda(\N)=
Ord(\N)$ is claimed in Theorem~\ref{Thy} and we prove it on page \pageref{proofy}. Theorem~\ref{Thb} is
contained in Theorem~\ref{ThIIIt} and  Corollary 2 of
Theorem~\ref{LemmaIIIoq}.

 We must
emphasize that in the case of proper Tarski hierarchy  the ordinal \emph{ $\lambda^p(\N)$ is much larger than the
first non-recursive ordinal in $\N$. } Consequently, the proper Tarski
hierarchy over $\N$ does not coincide with the sets hyperarithmetical
in $\N$. Though true, it is not therefore evident that the Tarski hierarchy does
stop at the first non-recursive ordinal. 

\section{General notions}
In this paper, We take \emph{ a structure }
to be a set of
predicates and function symbols where we assume predicates and function symbols
to be inherently interpreted. In addition, we assume predicates and function symbols to be
interpreted on the natural numbers $\omega$, i.e. the standard model
of natural numbers. Finally, we shall
deal only with structures of basic strength, i.e. those in which all
the usual arithmetical operations are definable.

\Definition{Defa}
\Benum
  \Item Let $n>0$. Then $P=\langle A,n\rangle$ is \emph{(n-ary) predicate }
iff $ A\subseteq \omega^n$;
$n$ will be called \emph{ the arity of P} and $A$ its \emph{ extension }.
$\langle A,n\rangle$ will also be denoted by $~^nA$.

  \Item Let $ n\geq0$. Then $ F=\langle f,n\rangle $ is  \emph{ (n-ary) function symbol } iff
$f:\omega^n\rightarrow \omega$ is a total $n$-ary function from $\omega^n$ to $\omega$;
if $n=0$ we assume $f\in \omega$;
$n$ will be called \emph{ the arity of $F$ } and $n$ its \emph{
  extension }. $\langle f,n\rangle$ will also be denoted by $~^nf$.

\Eenum
\Enddefinition

\Definition{Defb}  
\Benum
\Item \emph{The arithmetic},   $\mathbf A$, is the set
of predicates and function symbols, $\left\{=,<,S,+,.\right\}$,
interpreted in the usual way over $\omega$. 
\Item   \emph{M is a structure } iff
$\M$  is a set of predicates and function symbols
 and $\A \subseteq \M$
\Eenum
\Enddefinition

$\cal P,~ \cal F$ will denote the set of all predicate resp. function
symbols.
For a structure $\M$, we set $\cal P(M):=P\cap M$ and $\cal
P(F):=F\cap M$.
${\cal P}_n,~{\cal F}_n$ denotes the set of $n$-ary predicates
resp. of function symbols. ${\cal P}_n(\M),~{\cal F}_n(\M)$ is defined in
a similar fashion.  

If $Y$ is a set of predicates and function symbols then $\M+ Y$ denotes the structure
$\M\cup Y$.
If $Y=\left\{A_1\dots A_s\right\}$, we
shall write simply $\M+ A_1\dots A_s$

\emph{The first order variables}(or simply just \emph{variables}) 
are the elements of the set $\{x_i;i\in\omega\}$. The elements of the set
$\{X^k_i;i,k\in\omega\}$ are the \emph{second-order
 variables}, $k$ is the \emph{arity} of the variable $X^k_i$. For
binary \emph{logical connectives} we shall take
$\land,\lor,\equiv,\Rightarrow$ and $\lnot$ is the unary
connective. The symbols for quantifiers are $\exists$ and $\forall$.

Syntactical concepts, \emph{terms},
\emph{formulae} etc are defined in the usual way. 
\emph{Formula scheme} is simply a
second-order formula with no second-order quantifications; we shall
never need formulae of higher order.
A formula scheme $\psi$ will
be written as 
\[
\psi=\psi[Y_1,\dots Y_k](y_1,\dots y_n)=\psi[Y_1,\dots
Y_k]=\psi(y_1,\dots y_n)
,\]
where
$Y_1,\dots Y_k$ are the second and $y_1,\dots y_n$ are the first-order
variables occuring in $\psi$.  If $H_1,\dots
H_k$ are second-order variables or predicates of arities corresponding
to $Y_i$`s, then $\psi[H_1,\dots H_k]$ denotes the result of
substituting $H_i$ for $Y_i$ in $\psi$, $i=1,\dots k$. We may also
write only $\psi[Y_s/H_s]$.In an obvious way, we define the  $\psi[Y_s/\phi]$ 
where $\phi$ has $n$ free variables and the  arity of $Y_s$ is $\leq n$\footnote{ Here,
we must make sure that $\phi$ is substitutable in $\psi$, i.e. there
is no confusion between variables in $\psi$ and $\phi$};
 the extra variables in $\phi$ will serve as parameters.

 The class of all formulae
with $n$ free variables (resp. the class of formulae with $n$ free
variables of a structure $\M$) will be denoted by 
 $Fle_n$ (resp. $Fle_n(\M)$) 
If $\psi=\psi[Y_1,\dots Y_k](y_1,\dots y_m) $ where $Y_i$ is
of arity $n_i$, $i=1,\dots k$, then we write
$\psi\in Fle^{n_1,\dots n_k}_m$.
or alternatively $\psi\in Fle^{n_1,\dots n_k}_m(\M)$. 

\Definition{def2x}
Let $\mu\in\omega^{<\omega}$, $\mu= \langle k_1,\dots
k_s\rangle$.
\Benum
\Item we say that $\psi\in Fle^\mu_k$ iff $\psi\in Fle^{k_1,\dots
        k_s}_m$.
\Item   we say that $X\in
 P(\omega^\mu)$ iff $X=X_1,\dots X_s$ and for every $i=1,\dots s$,
   $X_i\subseteq \omega^{k_i}$. 

\Item If $X\in
 P(\omega^\mu)$ then 
$~^\mu X$ will denote the list of predicates $~^{k_1}X_1,\dots
~^{k_s}X_s$. 
\Eenum
\Enddefinition

The definitions of a formula being \emph{ true } or \emph{ satisfied }
by a sequence of natural numbers will be left to the
reader. \footnote{Note that formulae are taken interpreted in
  themselves. Hence  we do not say that $\psi$ is true in a
  structure $\M$, but simply that $\psi $ is true.} 

In the obvious manner we introduce partial function 
\[
Val:Term\times \omega^{<\omega}\rightarrow\omega
\]
such that if $t=t(x_{i_1},\dots x_{i_k})$, $i_1<\dots i_k$, then
$Val(t,\langle a_1,\dots a_k\rangle)=a$ iff $t(\overline{a_1},\dots
\overline{a_k})=\overline{a}$ is true.
\footnote{$\overline n$ denotes the $n$-th numeral}.
 
\Definition{Defe}
Let $\M,\cal N$ be structures.
\Benum\Item Let $\psi=\psi(x_{i_1},\dots x_{i_n})\in Fle^n$,
$i_1<i_2<\dots i_n$. Then \[
Ext(\psi)=\{\langle a_1,\dots a_n\rangle\in
\omega^n;\langle a_1,\dots a_n\rangle~satisfies ~\psi\}
\]
\Item We say that $X\subseteq \omega^k$ \emph{is defined by $\psi\in
        Fle^k$ } iff  $X=Ext(\psi)$. $X$ \emph
      {is definable in $\M$ } iff there is $\psi\in Fle^k(\M)$ which
      defines $X$. 

\Item The set of all $X\subseteq \omega$ definable in $\M$ will be
      denoted by $\cal D(\M)$.
\Item We say $\M\sim \cal N$ iff
      $\D(\M)=\D(\N)$.
      The classes of equivalence of the relation $\sim$ will be called
      \emph{ definability classes.}
\Item Let $F:P(\omega^{k_1,\dots
        k_s})\rightarrow P(\omega^s)$. Then $F$
      \emph{ is definable in $\M$} iff there is $\psi\in
      Fle^{k_1,\dots k_s}_n(\M)$ such that for every $X_1,\dots X_s\in
      P(\omega^{k_1,\dots k_s})$ and $X\in P(\omega^n)$ there is
      $F(x_1,\dots X_s)=X$ iff $X=Ext(\psi[~^{k_1}X_1,\dots
        ~^{k_s}X_s])$.

\Eenum
\Enddefinition
Since we assume that structures have at least the strength of
arithmetic we can find a
  simple coding function 
\label{code}
 \[
[~]:\omega^{<\omega}\rightarrow \omega
\]
which enables us to express quantification over finite sets and
sequences of numbers.
For a sequence $a_1,\dots a_n$ the number $[a_1,\dots a_n]$ will be
called \emph{ the code } or \emph{ the G\"odel number } of the
sequence $a_1,\dots a_n$. 
For $A=\{a_1,\dots a_n\}$, $[A]$ will denote the code of the sequence
$b_1,\dots b_n$ such that $A=\{b_1,\dots b_n\}$ and $b_1<b_2<\dots
<b_n$. If $S=s_1,\dots s_n$, where $s_i$ are sequences or finite sets
of numbers then $[S]:=[[s_1],\dots [s_n]]$.

An important
consequence is that inductively specified sets are definable, as we
state in the following lemma.

\Lemma{LemmaIA} Let $k>0$.
There exists $S^{k}$ a $\Sigma_0$-definable function in $\A$,
$S^k:P(\omega)\times (P(\omega^2)^k)\rightarrow P(\omega)$ with the following property:
\noindent
Let $C\subseteq \omega$. Let $R_1,\dots R_k\subseteq \omega^2$ be a list of binary relations.
Then for every $a\in \omega$
$a\in S^k(C,R_1,\dots R_k)$
iff $a$ is the code of a sequence $y_0,\dots y_n\in \omega$ such that
for every $j\leq n$ either
 \Benum
\Item $y_j\in C$, or
         \Item
 there is $1\leq l\leq k$
and $i_1,\dots i_s<j$ such that $R_l([y_{i_1},\dots y_{i_s}], y_j)$
\Eenum
\EndLemma
\Proof Easy. \Endproof
\section{Truth and universal predicates} 

We have introduced notions which describe semantics and
syntax of a structure. The notions are set-theoretical and hence
they cannot be directly taken as predicates or functions which
are assumed to range over natural numbers.  In order
to be able to define something like `the jump operator` we must
formulate concepts which describe properties of a structure
by means of predicates defined on natural numbers. 
For this purpose we define (proper) universal predicate
for $\M$ and the truth predicate for $\M$ under a coding $c$, $Tr_{\M,c}$.

 For a relation $R\subseteq \omega^2$ and $a\in \omega$, $R(a,.)$
will denote the set $\left\{x\in \omega;R(a,x)\right\}$. For relations of
bigger arity similarly.

\Definition{defIbm} \m. $P,G\subseteq \omega^2$.
\Benum
\Item  $G$ is \emph{ a universal set for $\M$ } iff for every $X\in {\cal
        D}(\M)$ there exists $n\in \omega$ such that $X=G(n,.)$, i.e., iff ${\cal
      D(M)}\subseteq \{G(n,.);n\in\omega\}$. $~^2G$ will be called
      \emph{a universal predicate for $\M$.}
\Item \emph{ P is a proper universal set for $\M$ }  iff $P$ is a universal set
      and for every $n\in \omega$ the set $P(n,.)$ is definable in $\M$, i.e. iff ${\cal
      D(M)}= \{G(n,.);n\in\omega\}$. $~^2P$ will be called
      \emph{a proper universal predicate for $\M$.}

\Item Let $G$ be a universal  set for $\M$,  let
      $X\subseteq \omega$. 
      Then $n\in \omega$ will be called \emph{ a $G$-code of X } iff
      $X=G(n,.)$. If $X\subseteq \omega^k$, $k>1$, then $n$ is \emph{
        a $G$-code of $X$ } iff $n$ is the $G$-code of the set
      $\{[a_1,\dots a_k];a_1,\dots a_k\in X\}\subseteq \omega$. 
\Eenum
\Enddefinition

We can view a
 universal set as a list of subsets of $\omega$
$G(0,.),G(1,.),G(2,.)\dots$ 
such that every definable set in $\M$
occurs in this list. If $G$ is a proper universal set then also
every member of that list is definable in $\M$. 
Consequently, a proper universal predicate enables us to express quantifications
over definable sets in $\M$, while the universal set enables us to express
quantifications over a class containing all definable sets in $\M$.

\Proposition{Thd}
\m. Let $G$ be a universal predicate for $\M$.  Then 
\Benum
\Item every set definable in $\M$
      is $\Delta_1$-definable in $\A+G$,
\Item $G$ is not definable in $\M$.
\Eenum
\EndProposition
\Proof
 1) is obvious. 2) is well-known. 
\Endproof

\bigskip
A (proper) universal predicate for $\M$ determines what are the
definable sets in $\M$, but does not show what is the internal
structure of $\M$, what predicates and functions are in $\M$ etc. 
On the other hand, the notion of truth
predicate for $\M$ under a coding $c$ which we introduce below is a complete description
of $\M$.  Two structures which define the same
sets, $\M_1\sim \M_2$, have the same (proper) universal predicates but in
general will possess different truth predicates. This relation between
truth predicate and proper universal predicate is expressed in the Proposition
\ref{ThId}.

\Definition{Defx}
 Let $\M$ be a structure.
  A one-to-one function $c:\M\rightarrow \omega$
          will be called \emph{ a  coding for $\M$.} 

Let $c$ be a coding for $\M$.
Then $\overline c:\M\cup(\mbox{ logical symbols })\rightarrow \omega$ is
the one-to-one function such that: 
i) if $x\in {{\cal P}_n}(\M)$ then $\overline c(x)= 2[c(x),n,0]$ ii) if
$x\in {{\cal F}_n}(\M)$ then $\overline c(x)= 2[c(x),n,1]$ iii) if $x=X^i_n$ is a second-order variable then 
                    $\overline c(x)= 2[i,n,2]$ iv)          
 $\overline c:\land, \lor,\lnot, \Rightarrow, \equiv,
                    \forall, \exists, (, )
                    \rightarrow 1,3,5,7,9,11,13,15$ respectively and 
              if $x=x_i$ is a first-order variable then 
                    $\overline c:x_i\rightarrow 2i+17$.
                             
\Enddefinition

Let $\M$ be a structure, $c$ a coding for $\M$. If $s_1,\dots s_n$ are
logical symbols or elements of $\M$ then
\[
[s_1,\dots s_n]_c\in \omega
\]
  will denote 
the number $[\overline{c}(s_1),\overline{c}(s_2),\dots \overline{c}( s_n)]$ and it
will be called \emph{ the c-G\"odel number of }, or simply \emph{ the $c$-code of
  $S_1,\dots s_n$}.

\Definition{defIbk}  \m. Let $c$ be a coding for $\M$.
\Benum 
\Item Let $X$ be a set of strings of symbols from $\M$ or logical
      symbols. Then $X_c:=\{[x]_c;x\in X\}\subseteq \omega$. 
\Item $Tr_{\M,c}\subseteq \omega$ is the set of $c$-G\"odel numbers of true
 sentences of $\M$. 
The predicate $^1Tr_{\M,c}$,
will be called \emph{ the truth predicate for $\M$ under the coding
  c.}
\Item $Tr^k_{\M,c}\subseteq \omega$ is the set of $c$-G\"odel numbers of true
 sentences of $\M$ which are in $\Pi_k$ or $\Sigma_k$ prenex form.
The predicate $^1Tr^k_{L,c}$,
will be called \emph{ the $k$-truth predicate for $\M$ under the coding
  c.} 
\Item $Val_{\M,c}\subseteq \omega^2$ is the relation such that
      $Val_{\M,c}(a,b)$ iff $a$ is a $c$-code of a closed term $t$ and
      $Val(t)=b$.
\Item $D_{\M,c}\subseteq \omega^2$ is the relation such that 
      \Benum
      \Item for every $X\in \M$
            $D_{\M,c}(\overline{c}(X),.)\subseteq \omega$ is the set of codes of $n$-tuples
            $\langle a_1,\dots a_n\rangle$ such that $\langle a_1,\dots
            a_n\rangle\in Ext(X)$
            (where the arity of $X$ is $n$ if $X\in {\cal P}_n$, and $n-1$ if $X\in
            {\cal F}_{n-1}$).
  
      \Item if $x\not\in Rng(\overline c\lceil \M)$ then $D_{\M,c}(x,.)=\left\{1\right\}$. 
      \Eenum      
\Eenum
\Enddefinition

The relation $D_{\M,c}$ determines what are the predicates and
functions of $\M$, what are their codes, arities and extensions.
 We may notice that in $\A+~^2D_{\M,c}$ we are
able to define the truth on the \it atomic \rm propositions in $\M$,
while the predicate $\Tr$ is not in general definable in $\A+~^2D_{\M,c}$, as we
shall see. 

\Proposition{Thc}
\m, $c$ be a coding for $\M$.
\Benum 
\Item The following are $\Delta_1$-definable in $\A+~^1Rng(\overline{c})$:
$Term(\M)_c$, $Fle(\M)_c$, $Fle_n(\M)_c$.
\Item $D_{\M,c}$ is $\Delta_1$-definable in $\A+\T$ and $Rng(\overline{c})$
      is  $\Delta_1$-definable in $\A+D_{\M,c}$
\Item $Val_{\M,c}$ is $\Delta_1$-definable in  $\A+~^2D_{\M,c}$.
\Item $Tr_{\M,c}^k$ is $\Delta_{k+1}$-definable in $\A+~^2D_{\M,c}$.
\Eenum
\EndProposition
\Proof 1), 3) and 4) are an easy application of Lemma~\ref{LemmaIA}.
2) is immediate.
\Endproof


\Proposition{ThId} \m~ and $c$ a coding for $\M$. Then
\Benum
\Item  there exists a proper universal  predicate for
      $\M$ which is $\Delta_1$-definable in $\A+\T$.

\Item Let $G$ be a universal predicate for $\M$.
      Then  $Tr_{\M,c}$ is definable in $\M+G+~^2D_{\M,c}$.
\Eenum
\EndProposition

\Proof 1) The relation $P$: 
\noindent
$P(n,x)$ iff $n=[\psi]_c$ ,
$\psi\in Fle_1(\M)$ and $\psi$ is in $\Sigma_k$ or $\Pi_k$ prenex form
 and $[\psi(\overline{x})]_c\in \Trk$ 
\noindent
is $\Delta_1$-definable and
it is a universal set for $\M$. 

2) The proof is an application of Lemma~\ref{LemmaIA} and proceeds as
follows. 

For a formula $\psi\in Fle(\M)$, a sequence $a_1,\dots
a_k\in\omega$ will be called a \emph{formula derivation for $\psi$}
iff i) $a_i$ is a $c$-code of a string $\psi_i$, $i=1,\dots i$, and
$a_k=[\psi]_c$ and 
ii) for every $i\leq k$ either $\psi_i$ is an atomic formula or there
are $i_1,i_2<i$ and $\psi_i=(\psi_{i_1})\triangle(\psi_{i_2})$, where
$\triangle$ is a binary logical connective, or $\psi=\triangle (\psi_{i_1})$,
where $\triangle$ is $\lnot$ or $\exists y$, $\forall y$.
  
A sequence $a_1,e_1\dots
a_k,e_k\in\omega$ will be called a \emph{ truth derivation for $\psi$ }
iff i) $a_1,\dots a_k$ is a formula derivation for $\psi$
and ii) if $a_i=[\psi_i]_c$, $\psi_i\in Fle_n$ then $e_i$ is a 
$G$-code of the set $\{[s_1,\dots s_l];l\geq n, \langle s_1,\dots s_n\rangle\in\omega^{<\omega}
\mbox{ satisfies }\psi\}$.

The proof of Proposition~\ref{Thc},1) requires to show that every formula of $\M$
has a formula derivation and the set of codes of
formula derivations is definable in $\A+~^2Rng(\overline{c})$. Here, it
must be shown that every formula of $\M$ has a truth derivation and
that the set of codes of truth derivations is definable in
$\M+G+^2D_{\M,c}$. Both parts are straightforward.
Finally, a $[\psi]_c\in \Tr$ iff $[\psi]_c\in Fle_0(\M)_c$ and $\psi$
has a truth derivation $a_1,e_1\dots
a_k,e_k$ such that $G(e_k,.)\not=\emptyset$.
\Endproof

\bigskip
\Corollary{~} Let $\M$ be a structure, $c$ a coding for
$\M$. Then 
\Benum
\Item Every set definable in $\M$ is $\Delta_1$-definable in $\A+\T$.
\Item $\Tr$ is not definable in $\M$.
\Eenum
\Endcorollary
\Proof Follows from the previous Proposition and Proposition~\ref{Thd}.
\Endproof

\Definition{defId} Let $\M$ be a structure. 
\Benum \item \emph{$\M$ is $L$-finite } iff $|\M|<\omega$,
i.e. iff $\M$ is a finite set function symbols and predicates. 

\item
$\M$ is \emph{
  essentially finite }  iff there
exists a structure $\M^\prime$ which is finite and $\M\sim \M^\prime$.

\Eenum
\Enddefinition

\bigskip
The following lemma expresses the key property of $L$-finite structures.

\Lemma{LemmaIk} Let $\M$ be a $L$-finite structure.
Then $D_{\M,c}$ is $\Delta_1$-definable in $\M$.
\EndLemma

\Proof  Let ${\cal P}(\M)=P_1,\dots P_s$, ${\cal F}(\M)=F_1,\dots F_k$.

For $P_i\in {\cal P}_n(\M),i\leq s$ there is $\psi_i$ a
$\Delta_1$-formula in $\M$ such
that for very $a\in \omega$, $a\in Ext(\psi_i)$ iff $a$ is a code of $n$-tuple
$a_1,\dots a_n$ and $\langle a_1,\dots a_n \rangle \in Ext(P)$.
Analogically, if $F_i\in {\cal F}_n(\M),i\leq k$ then there is $\psi_{s+i}$ a
$\Delta_1$-formula such that for every $a\in \omega$, $a\in Ext(\psi_{s+i})$
iff $a=[a_1,\dots a_{n+1}]$ and $\langle a_1,\dots a_n,a_{n+1} \rangle \in Ext(F)$.

Let $t_1,\dots t_s$ and $t_{s+1},\dots t_{s+k}$ denote the numerals corresponding to
$\overline{c}(P_1),\dots$ $\overline{c}(P_s)$ and
$\overline{c}(F_1),\dots$ $ \overline{c}(F_k)$.
Then $D_{\M,c}(x,y)$ is defined in $\M$ by the following formula

$((x\not=t_1)\land \dots (x\not=t_{s+l})\land y=\overline{1})
\lor ((x=t_1\land \psi_1(y)\lor \dots (x=t_{s+l}\land \psi_{k+l}(y)))$
\Endproof

\bigskip
 \Corollary{~}  Let $\M$ be a $L$-finite structure, $c$ a coding for $\M$, $k\in \omega$.
Then $Tr^k_{L,c}$ is $\Delta_{k+1}$-definable in $\M$. The sets 
$Term(\M)_c$, $Fle(\M)_c$, $Fle_n(\M)_c$ are $\Delta_1$-definable in $\M$.
\Endcorollary

\Proof Follows from the previous Lemma and Proposition~\ref{Thc}.\Endproof
\bigskip

For a given structure, by different choices of coding $c$ we
can obtain different truth predicates, and the structure $\M+\T$ will have
different expressive powers. Similarly for (proper) universal
predicates; in particular, if $\M$ is a structure and  $B\subseteq \omega$ is any
given set then we can find a (proper) universal predicate for $\M$
such that $B$ is definable in  $\M+~^2G$.  We see that neither the
universal nor the proper universal
predicate can have the role of `the jump operator` for $\M$, for such an
operation would not be unique. It is then an expectable move to try to choose a
particular (proper) universal predicate which would be
in some sense the weakest. This is achived using the concepts of
\emph{canonic universal predicate } and \emph{canonic proper universal predicate } which have been defined
on page \pageref{canonic}.\footnote{\label{defIc}
Note that we do
not introduce the symmetric concept of \it canonic truth predicate. \rm
The reason is that if we defined the Tarski hieararchy (see 
page \pageref{Tarski}) using the canonic truth predicate then 
the Theorem~\ref{ThIIIah} is false, ie. there would exist many
incomparable hierarchies over $\N$. In particular, for any $B\subseteq
\omega$ we could find a Tarski hierarchy $\{M_\alpha\}_{\alpha<\lambda(\N)}$ (defined in terms of canonic
truth predicate) such that $\omega<\lambda(\N)$ and $B$ is definable
in $\M_\omega$.}

\Lemma{LemmaId} \m.
\Benum
\Item  Assume that there is a coding $c_0$ for $\M$
       such that for every (proper) universal predicate $P$ the set
       $D_{\M,c_0}$ is definable in $\M+P$.
       Then $\M$ has a canonic (proper) universal predicate and if
       $P_0$ is a canonic (proper) universal predicate then
       $\A+P_0\sim \A+~^1Tr_{\M,c_0}$. 
            
\Item   Let $\M^\prime$ be a structure such that
$\M\sim \M^\prime$. Then $\M$ has a canonic (proper) universal predicate iff
$\M^\prime$ has a canonic (proper) universal predicate. If $P$ and $P^\prime$
are canonic (proper) universal predicates for $\M$ and $\M^\prime$ respectively
then $\A+P\sim \A+P^\prime$. 
\Eenum
\EndLemma
\Proof 1) follows from Proposition  \ref{ThId}. 2) follows from the fact
that $\M$ and $\M^\prime$ have the same (proper) universal predicates.
\Endproof

\Proposition{ThIk} Let $\M$ be an essentially $L$-finite structure. 
Then $\M$ has both a canonic and a canonic proper universal
predicate. If $P_0$ is a canonic (proper) universal predicate and $\N$ is a $L$-finite structure such that $\M\sim\N$ and $c$
a coding for $N$ then $\A+P_0\sim \A+~^1Tr_{\N,c}$.
\EndProposition
\Proof  By  Lemma~\ref{LemmaId} it is sufficient to show that $D_{\N,c}$ is definable
in $\M$. But that is claimed in Lemma~\ref{LemmaIk}. 
\Endproof 

\bigskip
Recall the definitions of Tarski and proper Tarski hierarchy given on
page \pageref{Tarski}.
Since for a given structure there in general exist infinitely many
canonic (proper) universal predicates, neither the Tarski
hierarchy nor the proper Tarski hierarchy are defined uniquely. 
The following Theorem shows that the hierarchies are unique at least
up to the equivalence $\sim$. 

 \Theorem{ThIIIah}  Let  $\cal N$ be a $L$-finite structure. Let $\{
{\cal M}_\alpha\}_{\alpha\in \xi_1}$ ,  $\{
{\cal M^\prime}_\alpha\}_{\alpha\in \xi_1}$ be two Tarski hierarchies over
$\cal N$.
Then $\xi_1=\xi_2>0$ and for every  $\alpha\in \xi_1$ there is
$\M_\alpha\sim\M^\prime_\alpha$. The same is true for two proper
tarski hierarchies.
\EndTheorem
\Proof Since $\J_0=\M_0=\N$ are finite then $\J_0,\M_0$ have canonic
proper universal predicates (Corollary of Proposition~\ref{ThIk}) and therefore $\xi_1,\xi_2>0$.
 The rest follows from Lemma~\ref{LemmaId},2).
\Endproof.

\Definition{defIIIt} Let $\cal N$ be a $L$-finite structure, $\{
\M_\alpha\}_{\alpha<\xi}$ be a Tarski hierarchy.
Then $\lambda(\J):=\xi$. If $\{
\M_\alpha^p\}_{\alpha<\xi}$ is a proper Tarski hierarchy
then $\lambda(\N)^p:=\xi$. 
\Enddefinition

A priori, we see that $\lambda(\N)$ and $\lambda^p(\N)$ can at
most be equal to $\aleph_1$, the first  uncountable ordinal. For then the structure $\M_{\aleph_1}$
is uncountable and there exist no truth or proper universal predicate for $\M_{\aleph_1}$
and we cannot hope to extend the hierarchies above $\aleph_1$.
The crucial question concerning the Tarski hierarchy and 
proper Tarski hierarchy is this:
\it is $\lambda(\J)$ countable\rm? If it is then the structure
$\M_{\aleph_1}$ is a countable structure which does not have a canonic
proper universal predicate and the proper Tarski hierarchy cannot be extended above
$\lambda(\J)$. If $\lambda(\J)=\aleph_1$ then we may say that the Tarski hierarchy
does not have an upper bound. 

\Theorem{ThIIIr} Let $\N$ be a $L$-finite structure. Let
$\{ \M_\alpha\}_{\alpha< \lambda(\J)}$ be a Tarski hierarchy over $\J$. Let
 $\alpha<\lambda(\J)$. Then
$\M_{\alpha+1}$ is essentially finite.
Hence
 $\alpha+1<\lambda(\J)$ and $\lambda(\J)$ is a limit ordinal.
The same is true for the proper Tarski hierarchy.
\EndTheorem
\Proof
$\M_{\alpha+1}=\M_\alpha+P$, where $P$ is a universal predicate.
But $\M_\alpha+P\sim_{\cal D}\A+P$, from Proposition~\ref{Thd}, 1). Hence $\M_{\alpha+1}$ is
essentially finite, it has a canonic universal predicate and
$\alpha+1<\lambda(\J)$
\Endproof

\section{Ordinals and the first part of Theorem~\ref{Tha}}
In this section we will prove that for a (proper) Tarski hierarchy
over $\N$ there is $Ord(\N)\leq \lambda(\N)$ (resp. $Ord(\N)\leq \lambda^p(\N)$).

\Definition{Def4a}
\m. Let $\lambda,\mu \in \omega^{<\omega}$.
\Benum
\Item $\cal A\it\subseteq P(\omega^\lambda)$ is \emph{ a system
defined by $\psi\in Fle^\lambda_0$ }  iff 
$$
{\cal A}=\{ A\in P(\omega^\lambda); A 
\mbox{ \it satisfies } \psi\}
$$
       $\cal A$ will be called \emph{ a definable system in $\M$ } iff it is
    defined by some $\psi\in Fle^\lambda_n (\M)$

\Item
Let  $A\in P(\omega^\lambda)$. Then $\psi\in Fle^\lambda_n$ is \emph{a proper implicit definition of A  } iff
$\psi$ defines the system $\left\{A\right\}$. 

\Item $B\in P(\omega^\mu)$ is \emph{ implicitly definable
        in $\M$ } iff there exist
$A\in P(\omega^\lambda)$, $A$ has a proper
implicit definition in $\M$ and $B$ is definable in the structure
$\M+~^\lambda A$.

\Item Let $F:P(\omega^\lambda)\rightarrow P(\omega^\mu)$,
      $\mu=\mu_1,\dots \mu_n$.
We will say that $F$ \emph{ is
defined by $\psi_1,\dots\psi_n$ }, $\psi_i\in Fle^\lambda_{\mu_i}$,  iff for every
$X=X_1,\dots X_s\in P(\omega^\lambda)$, $F(X)=\langle Y_1,\dots
Y_n\rangle$, we have
$Y_i=Ext(\psi_i(~^\lambda X)),i=1,\dots n$. That
$F$ is \emph{ definable in $\M$ } we introduce in the obvious way.

\Item Let $F: P(\omega^{\lambda})\rightarrow  P(\omega^{\mu})$.
We will say that $\psi\in
Fle^{\lambda,\mu}_0$ is \emph{a proper implicit definition of $F$ }
iff 
for every $X\in  P(\omega^{\lambda})$
there is a unique $Y\in  P(\omega^{\mu})$ such that $\psi[~^\lambda X,~^\mu Y]$ is
true and for such $Y$, $F(X)=Y$.
                 
\Item Let $F:P(\omega^\lambda)\rightarrow P(\omega^\mu)$.
We will say that \emph{ $F$ is implicitly definable in $\M$ } iff there are
functions $F_1,F_2$, $F_1:P(\omega^\lambda)\rightarrow P(\omega^\pi)$,
$F_2:P(\omega^\pi)\rightarrow P(\omega^\mu)$ such that $F(X)=F_2(F_1(X)),~X\in N^\lambda$
and $F_1$ has a proper implicit definition in $\M$ and $F_2$ is
definable in $\M$.
\Eenum
\Enddefinition

We may observe that

\Benum
\Item
If $B\in P(\omega^\lambda)$,  $F:P(\omega^\lambda)\rightarrow
P(\omega^\mu)$ 
are definable in
$\M$ then they have a proper implicit definition in $\M$.
If they have a proper implicit definition in $\M$ then they are
implicitly definable in $\M$.
\Item Let $F_1:P(\omega^\lambda)\rightarrow P(\omega^\pi)$,
$F_2:P(\omega^\pi)\rightarrow P(\omega^\mu)$ and $F(X)=F_2(F_1(X))$, 
$X\in \omega^\lambda$.
Then
\Benum
\Item if $F_1,F_2$ are definable resp. implicitly definable in $\M$
then $F$ is definable resp. implicitly definable in $\M$.
         
\Item if $F_1$ is definable in $\M$ and $F_2$ has a proper implicit
definition in $\M$ then $F$ has a proper implicit definition in $\M$
\Eenum
\Item if $B\in P(\omega^\lambda)$ and $F:P(\omega^\lambda)\rightarrow P(\omega^\mu)$ are
definable resp. implicitly definable in $\M$ then $F(B)$ is definable
resp. implicitly definable in $\M$.

\Eenum

The following statement will not be used in this work but it gives
an important characterisation of implicitly definable sets. We
therefore do not enter the proof.

\bigskip
\noindent
\bf Proposition. \it \m. Then $B\in P(\omega^k)$ is implicitly
definable in 
$M$ iff it is $\Delta^1_1$ in $\M$ (i.e. iff $B$ is hyperarithmetical in
$\M$). \rm
\Proof
The implication '$\rightarrow$' is obvious. The other follows from
Lemma \ref{ThIIaa}. 
\Endproof

\Lemma{LemmaIIa} \m, $\lambda,\mu\in \omega^{<\omega}$.
\Benum
\Item Let $B\in P(\omega^\lambda)$ be implicitly definable in $\M$. Then there exists
$A\subseteq \omega$, $A$ has a proper implicit definition in $\M$ and $B$ is
definable in $\M+~^1A$
\Item Let $B\in P(\omega^\lambda)$ be implicitly definable in $\M$ ,
      $C\in P(\omega^\mu)$ 
implicitly
definable in  $\M+ B$. Then $C$ is implicitly definable in $\M$.
\Eenum
\EndLemma
\Proof Straightforward.
\Endproof

\Definition{defIIa} Let $\M$ and $\N$ be structures. Then
\Benum
\Item
$\M$ is \emph{ implicitly closed } iff every set which is implicitly
definable in $\M$ is definable in $\M$.
\Item
 $\cal I(M)$ is the structure $\M+\{~^1X;
X\subseteq \omega,~ X \mbox{ \it implicitly definable in } \M\}$. 
\Eenum
\Enddefinition

\Corollary{of Lemma~\ref{LemmaIIa} }
\m. Let $\N:=\cal I(\M)$. Then  
i) $\N$ is implicitly closed, ii) $\cal D(\M)\subseteq\cal D(\N)$ and iii) for every
$\N^\prime$ if $\N^\prime$ satisfies i) and ii) then  $\D(\N)\subseteq \D(\N^\prime)$.
\Endcorollary

\Proof
Let $\M$ be given.  By Lemma~\ref{LemmaIIa}, 2) if a set is
implicitly definable in $\cal I(M)$ then it is implicitly definable in $\M$.
Hence $\cal I(M)$ is implicitly closed. The rest is immediate.
\Endproof

\Proposition{ThIIa}  There is a function $TR:P(\omega^2)\rightarrow
P(\omega) $ which has a proper implicit definition in $\A$ such that
for every structure $\M$ and a coding $c$ for $\M$ we have 
\[
TR(D_{\M,c})=\Tr
\]
\EndProposition
\Proof First, observe that Proposition~\ref{Thc},3) can be strengthened to
assert that there exists a function $VAL:P(\omega^2)\rightarrow
P(\omega^2)$ definable in $\A$ such that for every structure $\M$ and
a coding $c$ for $\M$, 
\[
VAL(D_{\M,c})=Val_{\M,c}
\] 
 Let $\M$ be a structure and $c$ its coding.
Then $X=\Tr$ iff $X\subseteq Fle_0(\M)_c$ and for every $x=[\psi]_c\in
Fle_0(\M)_c$ the following conditions are satisfied 
\Benum
\Item If $\psi=P(\overline{n_1},\dots \overline{n_i})$ is atomic and
      $\overline{n_l}$ is the $n_l$-th numeral then 
      $x\in X$ iff $D_{\M,c}(\overline{c}(P),[n_1,\dots n_l])$

\Item If $\psi=P({t_1},\dots {t_i})$ is a closed atomic formula, where
      $t_1,\dots t_i$ are terms, then 
      $x\in X$ iff $[P(\overline{Val_{\M,c}(t_1)},\dots \overline{Val_{\M,c}(t_i)}]_c\in X$.
\Item If  $\psi=\lnot \xi$ then $x\in X$  iff not $[\xi]_c\in X$.
If  $\psi=\xi_1\land\xi_2$ then $x\in X$ iff $[\xi_1]_c\in X$ and
$[\xi_2]_c\in X$ and so on for the other logical connectives.

\Item if $\psi=\exists y\eta$ then $x\in X$ iff there exists $a\in
\omega$ such that $[\eta(y/\overline a)]_c\in X$.
If $\psi=\forall y\eta$ then $x\in X$ iff for every $a\in
\omega$,  $[\eta(y/\overline a)]_c\in X$.
\Eenum
Let $S\in Fle^{2,1}(\A)$ be a formula scheme obtained as a natural
translation of the above conditions and by replacing every occurence
of $^2D_{\M,c}$ (including the one in $Val_{\M,c}=VAL(D_{\M,c})$) by
a second-order variable $Y$. Then we can see that $S$ is a proper
implicit definition of a function $TR$ with the desired property. 
\Endproof

\bigskip
\Corollary {1. }
Let $\M$ be a $L$-finite structure, $c$ a coding for $\M$.
Then the truth predicate $\Tr$ has a proper implicit definition in
$\M$.
\Endcorollary

\Proof It must be shown that $D_{\M,c}$ is definable in $\M$ if $\M$ is
finite. But that has been claimed in Lemma~\ref{LemmaIk}.
\Endproof

\bigskip
\Corollary {2. } Let $\M$ be an essentially $L$-finite structure.
Then there is a proper universal predicate for $\M$ which is implicitly definable in
$\M$. Hence, $\M$ is not implicitly closed.
\Endcorollary

\Proof Apply  Proposition~\ref{ThIk} on $L$-finite structure $\J$ such that
$\N\sim\M$ to show that a universal predicate for $\M$ is implicitly
definable in $\M$. That a universal predicate is not definable in $\M$
is claimed in Proposition~\ref{ThId}.
\Endproof

\bigskip
We shall see that one of the important characteristics of a structure is
how many ordinals are definable in the structure. 
We shall say that $R$ is a \emph{ linear ordering on $X$ } iff $R$ is
reflexive, transitive, and weakly antisymmetric on $X$ and for every $x,y\in
X$,  $R(x,y)$ or $R(y,x)$. $R$ is \emph{ a linear ordering }
iff $R$ is a linear ordering on $Rng(R)$. Thus we take a
linear ordering to be  \it non-strict. \rm In order to avoid
confusion, we shall also  write
$\preceq_R$ instead of $R$. $x\prec_R y$ is then defined as $x\preceq_R y$ and
$x\not=y$.
Note that for a linear ordering we have 
$Rng(R)=Dom(R)$.  
If $X\subseteq Rng(R)$ then we define $R\lceil X:=R\cap X^2$. 
If $n\in Rng(R)$ then $R_n$ will denote the relation such that
\begin{center}  $R_n(x,y)$ iff $R(x,y)$ and $R(y,n)$ \end{center}
\label{Rn}

\Definition{defIIe}  Let $\ro\subseteq \omega^2$ be a linear ordering, let $\alpha$ be a
countable ordinal.
\Benum
\Item Then $\ro$ is \emph{ a representation of ordinal $\alpha$ } 
      iff $\prec_\ro$ is a well-ordering of the order-type $\alpha$. 
\Item       Let $\beta< \alpha$. Then  $\ro_\beta$ will be defined by
      induction as follows:
      let $\ro_0:=a$, $a$ the $\ro$-smallest member of $Rng(R)$.
      If $\beta>0$, let $\ro_\beta$ be the $\ro$-smallest member of the set
      $Rng(\ro)\setminus\left\{\ro_\gamma;\gamma<\beta\right\}$.
\Item
      Let $\beta\leq \alpha$. Then $\ro\lceil\beta$ is the representation of $\beta$
      such  that  ${\ro\lceil\beta}=\ro\lceil
      \{\ro_\gamma;\gamma<\beta\}$.
\Eenum
\Enddefinition

Thus $\ro\lceil\beta$ is a representation of $\beta$.
$\ro_\beta$ is the $\ro$-smallest element majorising $Rng(\ro\lceil \beta)$ if some such $\ro_\beta$ exists (if $\beta=\alpha$ then
$\ro\lceil\alpha=\ro$ while $\ro_\alpha$ is not defined).

\Definition{defIIg} \m,
\Benum
\Item  Let $\alpha$ be a countable ordinal. Then
      $\alpha$ is \emph{ (implicitly) definable
        in $\M$ } iff there is $\ro$ a representation of $\alpha$ which
      is (implicitly) definable in $\M$.
\Item The smallest undefinable ordinal in $\M$ will be denoted by $Ord(\M)$.
\Eenum
\Enddefinition

We can see the following:
\Benum
\Item every $\alpha>Ord(\M)$ is undefinable in $\M$. I.e.,
      the set of definable ordinals in $\M$ is an interval.
\Item $0,1,\dots\omega$ are definable in $\M$.
\Item If $\alpha,\beta$ are definable in $\M$ then $\alpha+\beta$ and
      $\alpha.\beta$ are definable in $\M$. Hence $Ord(\M)$ is a limit
      ordinal.
\Eenum

Now we shall define two important concepts: the concept of iterated
truth predicate over a well-ordering, $\Trro$, and the  notion
of iteration of a general operation over a well-ordering.

\Definition{defIIh} 
Let $n>0$, $F:P(\omega^{n+1})\rightarrow P(\omega^n)$, let $B\subseteq N^n$.
Let
$\ro$ be a representation of an
ordinal $\alpha$. Let $Z\subseteq \omega^{n+1}$. For $\beta< \alpha$  we define
$$
Z_{<\beta}:=\{ \langle a_0,\dots a_n\rangle\in \omega^{n+1}; Z(a_0,\dots
a_n)~and~a_0\prec_{\ro}\ro_\beta \}
$$
We will say that $Z =REK_n(B,F,\ro)$ iff $Z$ satisfies
the following conditions
\Benum
\Item  If $y\not\in Rng(\ro)$ then $Z(y,.)=\emptyset$.  
If $\ro\not=\emptyset$, let $Z(\ro_0,.)=B$
\Item If $0<\beta<\alpha$ then $Z(\ro_\beta,.)=F(Z_{<\beta})$
\Eenum
\Enddefinition

We note that
 \Benum
\Item $Z=REK(B,F,\ro)$ as defined above exists
and is unique, 
    
\Item the definition of $REK_n(B,F,\ro)$ can be rewritten as a formula
scheme, as we state in the following proposition.
\Eenum

\Proposition{ThIId} Let $F:P(\omega^{n+1})\rightarrow P(\omega^n)$ have a proper implicit
definition in a structure $\M$.
Then there exists a function $REK_F:P(\omega^{n,2})\rightarrow
P(\omega^{n+1})$ which has a proper implicit definition in $\M$ with
the following property:  
\noindent
for every $B\subseteq N^n$ and $\ro$ a representation of an ordinal
\[
REK_F(B,\ro)=REK_n(B,F,\ro)
\]
\EndProposition

\Corollary{~}  \m, let $\ro,B\subseteq \omega^n,F:P(\omega^n)\rightarrow
P(\omega^n)$ let $Z:=REK_n(B,F,\ro)$. Then
if $B,F,\ro$  have a (proper) implicit definition in $\M$ then $Z$
has a (proper) implicit definition in $\M$.
\Endcorollary

\Proof Straightforward.
\Endproof

\bigskip
Later, we shall see that every
set which is implicitly definable in $\M$ is also definable in terms
of some $REK_n(B,F,\ro)$, where all $B,F,\ro$ are definable in $\M$.

\Definition{Def4a2}
Let $R$ be a linear ordering. \m, and $c$ a coding for $\M$. 
\Benum
\Item We will say that \emph{$R$ and $c$ are compatible } iff
      $Dom(c)\cap Rng(R)=\emptyset$
\Item For $U\subseteq \omega^2$ we shall write that $U\in      \overline{Tr}(\M,c,R)$ 
iff the following is satisfied: \Benum
   \Item   If $x\not\in Rng(R)$ then $U(x,.)=\emptyset$.
         If there is $R_0$ the $R$-first element of $Rng(R)$, 
then $U(R_0,.)=Tr_{\M,c}$
                 
   \Item If $R_0<_R n$ then
$U(n,.)=Tr_{\M_{<n},{c}_n}$,
where ${\M_{<n}}$ is the structure
$\M+\{~^1U(m,.);m<_Rn\}$
and $c_n$ is the coding induced on ${\M}_{<n}$ (by $U$
and $c$).
   \Eenum 
   \Item If $R$ is a well-ordering then $\Trro$ is the set such that 
$ \overline{Tr}(\M,c,R)=\{Tr_{\M,c,R}\}$.
\Eenum
\Enddefinition
Clearly, if $\ro$ is a well-ordering compatible with $c$ then 
$\Trro$ can be  defined  as an iteration of adding a truth predicate
along the well-ordering
$\ro$. In this case we have 
$\overline{Tr}(\M,c,\ro)=\{\Trro\}$.
We will see in the last section that $\overline{Tr}(\M,c,\ro)$ is non-empty even for linear
orderings which are not well-orderings; in that case $U$ will not in
general be unique. Here, we shall deal with
$\overline{Tr}(\M,c,\ro)$ only in the case $\ro$ is a well-ordering.  
The main results about $\Trro$ presented below are that i) it is strong enough to
define all sets of the form $REK(B,F,\ro)$, for $B,F$ being definable
and ii) we can characterise the Tarski hierarchy by sets of the form
$\Trro$ with $\ro$ definable in $\M$.

\Proposition{ThIIae} \Benum
\Item
There is $TR^\star:P(\omega^2)\rightarrow P(\omega)$
which has a proper implicit definition in $\A$ with the following property:
\noindent
let $\M$ be a structure,  $\ro$ a representation of
an ordinal $\alpha$ and 
$c$ a coding for $\M$ compatible with $\ro$.
Then $\Trro=REK_n(Tr_{\M,c},TR^\star,\ro)$.

\Item Moreover, there exists a function
      $TRO:P(\omega^{2,2})\rightarrow P(\omega^2)$ with a proper
      implicit definition in $\A$ with the following property:
\noindent
let $\M$ be a structure and $c$ a coding for $\M$. Let $\ro$ be a
representation of an ordinal such that $\ro$ and $c$ are compatible.
Then \[TRO(D_{\M,c},\ro)=\Trro\]
\Eenum
\EndProposition

\Proof For 1), use  Proposition~\ref{ThIIa} 
and 2) immediately follows.
\Endproof

\bigskip
 \Corollary{}  \m, Let $c$ be a coding for $\M$ compatible with $\ro$,
$\ro$ being a representation of an ordinal.
Then $\Trro$ has a proper implicit definition in
$\A+~^2D_{\M,c}+~^2\ro$.
\Endcorollary

\Lemma{LemmaIIac}  \m, $c$ a coding for $\M$. Let $\ro$ be a representation
of an ordinal $\alpha$, $\ro$ and $c$ compatible. Then
\Benum
\Item for every  $\beta<\alpha$ we have
$$\A+~^2Tr_{\M,c,\ro\lceil(\beta+1)}\sim\A+^1Tr_{\M,c,\ro}(\ro_\beta,.)$$

\Item $\ro\lceil\beta$ is definable in $\A+~^2Tr_{\M,c,\ro\lceil\beta}$.
\Item  
If $\beta\leq \alpha$ let us define
\[
\M_{<\beta}:=\A+\{
~^1Tr_{\M,c,\ro}(\ro_\gamma,.);\gamma<\beta\}
\] 
Let $c_\beta$ be the coding for $\M_{<\beta}$ induced on $\M_{<\beta}$.
Assume that $\beta$ is a limit ordinal. Then
there is a universal predicate for $\M_{<\beta}$ definable in $A+~^2Tr_{\M,c,\ro\lceil\beta}$.
\Eenum
\EndLemma
\Proof 1) and 2) are straightforward. In 3) notice that 
every set definable in $\M_{<\beta}$ is $\Delta_1$-definable in
$A+~^2Tr_{\M,c,\ro\lceil\beta}$ and that $1$-truth predicate
$Tr_{M_{<\beta},c_\beta}^1$ is definable in $A+~^2Tr_{\M,c,\ro\lceil\beta}$.
\Endproof

\Lemma{LemmaIIab} \m, let $P$ be a universal
predicate for $\M$. Let $B\subseteq
N^n$.
 Let $F:P(\omega^{n+1})\rightarrow P(\omega^n)$ have a proper implicit definition in $\M+P$.
Assume that $\beta$ is a limit ordinal and that $\ro$
is a representation of $\alpha\geq\beta$ such that $\ro\lceil\beta$ is
definable in $\M+P$. Assume that  for every $\gamma<\beta$,
$REK_n(B,F,\ro\lceil{(\gamma+1)})$ is  definable in $\M$.
Then $REK_n(B,F,\ro\lceil\beta)$ is definable in $\M+P$.
\EndLemma
\Proof 
Let $\theta$ be a proper implicit definition of $F$ in $\M+P$. Let $\eta$
be a definition of $\ro\lceil\beta$ in $\M+P$. 
Let $R_\theta$ be a formula scheme in $\M+P$ which is a proper
implicit definition of the function $REK_F$ (see Proposition~\ref{ThIId}).
Let $\eta^\prime(x,y,z)$ be the
formula $\eta(x,y)\land \eta(y,z)$. Then for every
$a=\ro_\gamma,\gamma<\beta$, we have $Ext(\eta^\prime(x,y,\overline
a))=\ro\lceil(\gamma+1)$, and if $a\not \in O_{\ro\lceil\beta}$ then
$Ext(\eta^\prime(x,y,\overline a))=\emptyset$.
Since
$REK_n(B,F,\ro\lceil\gamma)(\ro_0,.)=B$ then $B$ is definable in $\M$.
Let $\xi$ be a definition of $B$ in $\M$. Let $S$ be the scheme

$$S[X^{n+1}](z):=REK_\theta[X^{n+1},\xi,\eta^\prime]
$$
Then for every $a=\ro_\gamma$ $Z\subseteq \omega^{n+1}$ satisfies $S(\overline a)$
iff $Z=REK_n[B,F,\ro\lceil(\gamma+1)]$, and if $a\not\in O_{\ro\lceil\beta}$ then
$S(\overline a)$ is satisfied by $\emptyset$ only. 
In $\M$ we can define the
relation $Q\subseteq N^2$ such that $Q(a,b)$ iff $b$ is a $P$-code of
a set $Z\subseteq \omega^n$ which
satisfies $S(\overline a)$. Because we assumed that
$REK_n(B,F,\ro\lceil{(\gamma+1})),\gamma<\beta$ is definable in $\M$ and $P$
is a universal predicate for $\M$ then

i) for every $a=\ro_{\gamma+1},\gamma<\beta$, $Q(a,.)\not=\emptyset$
and furthermore

ii) if $m\in Q(a,.)$ and $a=\ro_\gamma,\gamma<\beta$ then $m$ is a $P$-code of
$REK_n(B,F,\ro_\gamma)$, and if $a\not\in O_{\ro\lceil\beta}$ then
$P(m,.)=\emptyset$.

Hence the following are equivalent

a) $\langle k_1,\dots k_{n+1}\rangle \in REK_n(B,F,\ro\lceil\beta)$

b) there exist $a,b\in N$ $Q(a,b)$ and $[k_1,\dots k_{n+1}]\in P(b,.)$

But this equivalence can be written as a definition of
$REK_n(B,F,\ro\lceil\beta)$ in $\M+P$
\Endproof 

\Proposition{ThIIe}\m. Let $F:P(\omega^{n+1})\rightarrow P(\omega^n)$ and $B\subseteq
\omega^n$ be definable in $\M$. Let $Z=REK_n(B,F,\ro)$,
where $\ro$ is a representation of $\alpha$.
Let $c$ be a coding for $\M$ compatible with $\ro$.
Then $Z$ is definable in $\M+ \Trro$.
\EndProposition
\Proof
Let us prove by induction that for every $\beta\leq\alpha$,
$Z_\beta:=REK_n(B,F,\ro\lceil\beta)$ is definable in 
$\M_\beta:=\A+~^2Tr_{\M,c,\ro\lceil\beta}$.

Assume that $\alpha>0$, otherwise the proposition is trivial.

\noindent
We have
$Z_0=\emptyset$ and $Z_1=\{\ro_0\}\times B$ which are
definable in $\A$ and resp. in $\D(\M)\subseteq \D(\A+~^1\Tr)\sim (\A+~^2Tr_{\M,c,\ro\lceil
  1})$. 

Assume the statement holds for every $\gamma<\beta$.

Assume that $\beta$ is isolated. Then $Z_{\beta-1}$ is definable in
$\M_{\beta-1}$.
We have
\[
Z_\beta=Z_{\beta-1}\cup \{\ro_{\beta-1}\}\times F(Z_{\beta-1}) 
\]
But $F$ is definable in $\M$ and therefore $F(Z_{\beta-1})$ and hence $Z_\beta$ are definable
already in $\D(\M_{\beta-1})\subseteq \D(\M_\beta)$.

Assume that $\beta$ is a limit. 
By the assumption, every $Z_\gamma$, $\gamma<\beta$ is
definable in
$\N^\prime:=\A+\{~^2Tr_{\M,c,\ro\lceil\gamma};\gamma<\beta\}$. By Lemma
\ref{LemmaIIac}, 2) we have
$\N^\prime\sim \A+\{~^1Tr_{\M,c,\ro\lceil\beta}(\ro_\gamma,.);\gamma<\beta\}$
and hence every $Z_\gamma$, $\gamma<\beta$, is definable in $\M_{<\beta}$.
We shall apply Lemma~\ref{LemmaIIab}. Let us check that the
assumptions of the lemma are satisfied. 
By Lemma~\ref{LemmaIIac},1) $\ro\lceil\beta$ is definable in
$\M_\beta$. By Lemma~\ref{LemmaIIac},3) a universal predicate for $\M_{<\beta}$
is definable in $\M_\beta$. Hence, by Lemma~\ref{LemmaIIab}, $Z_\beta$
is definable in $\M_\beta$.
\Endproof

\Lemma{LemmaIIIos} Let $\J$ be a $L$-finite structure and $c$ a coding for $\J$.
Let $\ro$ be a representation of ordinal $\alpha$ compatible with $c$.
Let $\beta\leq\alpha$ and let
$\M_{<\beta}$, $c_\beta$ be as defined in Lemma
\ref{LemmaIIac},3). Let $P$ be a universal predicate for $M_{<\beta}$ such that
$\ro\lceil\beta$ is definable in $M_{<\beta}+P$.
Then $Tr_{M_{<\beta,c_\beta}}$ is definable in  $\M_{<\beta}+P$.
\EndLemma
\Proof 
Let us first show that $Tr_{\M,c,\ro\lceil\beta}$ is definable in $\A+P$.

 Assume that $\beta$ is isolated. Then
$~^1Tr_{\J,c,\ro}(\ro_{\beta-1},.)\in \M_{<\beta}$ and hence it is definable
in $\A+P$. But from Lemma {\ref{LemmaIIac},1) we have
$\A+~^1Tr_{J,c,\ro}(\ro_{\beta-1},.)\sim \A+~^2
Tr_{\M,c,\ro\lceil\beta}$ and $Tr_{\M,c,\ro}(\ro_{\beta},.)$ is
definable in $\A+P$.

Assume that $\beta$ is limit. We shall use Lemma~\ref{LemmaIIab} (note that
Proposition~\ref{ThIIae} asserts that $Tr_{\M,c,\ro}=REK(Tr_{\M,c},TR^\star,\ro)$ where
$TR^\star$ has a \emph{ proper } implicit definition).
From Lemma~\ref{LemmaIIac} we have

$$\M_{<\beta}\sim
\M+\{~^2Tr_{J,c,\ro\lceil(\gamma+1)};\gamma<\beta\}
$$
and hence every $Tr_{\M,c,\ro\lceil (\gamma+1)}$, $\gamma<\beta$, is definable in $\M_{<\beta}$.
Furthermore, since $P$ is a universal predicate for $\M_{<\beta}$, then
by
Lemma~\ref{LemmaIIab}, $Tr_{\M,c,\ro\lceil\beta}$ is definable in
$\M$.

It is trivial to show that the set $D_{\M_{<\beta},c_\beta}$ is
definable in $\A+~^2Tr_{\M,c,\ro\lceil\beta}$ and hence it is
definable in $\A+P$. Therefore, by Proposition~\ref{ThId},
$Tr_{\M_{<\beta},c_\beta}$
is definable in $\A+P$.
\Endproof

\Theorem{LemmaIIIoq} Let $\J$ be a $L$-finite structure. Let
$\{\M_\alpha\}_{\alpha< \lambda(\J)}$ be a Tarski hierarchy over $\J$. 
Let $\beta$ be a definable ordinal in $\J$.
\Benum
\Item Then $\beta<\lambda(\J)$
\Item Furthermore, let $c$ be a coding for $\J$. Let $\ro $ be a representation of
the ordinal $\beta+1$ compatible with $c$ definable in $\J$. Let $P_{\beta}$ be a
canonic universal predicate for $\M_{\beta}$. Then
$\A+P_{\beta}\sim\A+~^1\Trro(\ro_\beta,.)$. 
\Eenum
The same is true for the proper universal predicate and proper Tarski hierarchy. 
\EndTheorem
\Proof We shall say that $c$ is an \it ultracanonic coding
for a structure $\M$ \rm iff $\Tr$ is definable in every $\M+P$, where $P$ is a
universal predicate for $\M$. 
From Lemma~\ref{LemmaId} and Proposition~\ref{ThId} we obtain the following:
\it
\Smallskip
Let $\N_1\sim \N_2$. Assume that
$\N_1$ has an ultracanonic coding $c_1$. Then  $\N_2$ has a canonic
 universal predicate and if $P$ is a canonic universal predicate for
$\N_2$ then $\A+Tr_{\N_1,c_1}\sim \A+P$.
\Smallskip\rm

Assume that $\beta,\ro,\M,c$ are as in the statement 2).
By transfinite induction we shall prove the proposition:
\it\Smallskip
For every $\alpha\leq\beta$ it is the case that $\alpha<\lambda(\J)$. Moreover, if $P_\alpha$
denotes the canonic universal predicate for $\M_{\alpha}$
then $A+P_\alpha\sim\A+^1Tr_{\ro,\M,c}(\ro_\alpha,.)$
\rm\Smallskip
First, let $\alpha=0$. Then $\alpha<\lambda(J)$. From the definition of
$\Trro$ we obtain $Tr_{\N,c,\ro}(\ro_0,.)=\Tr$. Furthermore, since $\N$ is
finite then any coding for $\N$ is ultracanonic (Proposition~\ref{ThIk}). Hence
$\A+~^1\Tr\sim \A+P_0$ and so $\A+~^1Tr_{\J,c,\ro}(\ro_0,.)\sim\A+P_0$.

Let $0<\alpha\leq\beta$ and assume that the proposition is true for every
$\gamma<\alpha$. Let $\M^\prime$ be the structure
$\M_\alpha=\M+\{P_\gamma;\gamma<\alpha\}$. 
 By the assumption 
$\M^\prime\sim\M+\{~^1Tr_{\ro,\M,c}(\ro_{\gamma},.),\gamma<\alpha\}$. Let
$\M_{<\alpha}$ denote the structure on the right hand side and let
$c_\alpha$ be the coding induced 
on $\M_{<\alpha}$. 
By the
previous Lemma, $c_\alpha$ is an ultracanonic coding for $\M_\alpha$
and hence $M_\alpha$
has a canonic universal predicate and if  $P_\alpha$ is a canonic
universal predicate for $\M_\alpha$ then $\A+
Tr_{\M_{<\alpha},c_\alpha}\sim \A+P_\alpha $. But from the definition
of $\Trro$ we have
$Tr_{\M_{<\alpha},c_\alpha}= Tr_{\N,c,\ro}(\ro_\alpha,.)$; hence 
$\A+P_\alpha\sim\A+~^1 Tr_{\N,c,\ro}(\ro_\alpha,.)$.

For the proper Tarski hierarchy the proof is exactly the same.
\Endproof

\bigskip
\Corollary{1} Let $\{ \M_\alpha\}_{\alpha<\lambda(J)}$ be a
Tarski hierarchy over $\J$. Let $\beta$ be isolated,
let $c$ be a coding for $\M$. Let $\ro $ be a definable representation of
the ordinal $\beta$ in $\J$ compatible with $c$. Then
$$\M_{\beta}\sim\A+~^1Tr_{\J,c,\ro}(\ro_{\beta-1},.)
\sim\A+~^2Tr_{\J,c,\ro\lceil \beta}$$
\Endcorollary
\Proof Follows from the previous Theorem and Lemma~\ref{LemmaIIac}.
\Endproof

\bigskip
\Corollary{2} Let $\{ \M_\alpha\}_{\alpha<\lambda(\J)}$ be a
Tarski hierarchy over $\J$ and let \linebreak$\{ \M_\alpha^p\}_{\alpha<\lambda^p(\J)}$ be a
proper Tarski hierarchy over $\J$. Then for every $\beta\leq Ord(\J)$,
we have
$\beta\leq \lambda(\J)\cap\lambda^p(\J)$ and  $\M_\beta\sim \M_\beta^p$.
\Endcorollary

\Proof Immediate.
\Endproof

\section{ Trees and the second half of Theorem~\ref{Tha} }

We now proceed to prove the rest of Theorem~\ref{Tha}, i.e., to show that
$\lambda(\N)=Ord(\N)$ and that $\M_{\lambda(\N)}=\cal I(\N)$. We shall
first prove the theorem (see page \pageref{proofx})  
\Theorem{Thx} Let $\N$ be a $L$-finite structure,
let $\{M_\alpha\}_{\alpha<\lambda(\N)}$ be a Tarski hierarchy over
$\N$. Then $\M_{Ord(\N)}\sim{\cal I(N)}$.
\EndTheorem
 Second, we will
prove (see page \pageref{proofy})
\Theorem{Thy} Let $\N$ be a $L$-finite structure. Then  $\cal I(N)$ does not
have a canonic universal predicate.
\EndTheorem
 For those purposes, we
shall use some properties of trees and linear orderings.  Trees are a standard tool for proving uniformization
results (see for example [4]). Results concerning the
definability of well-orderings can be found in [2].

\Definition{defIIad} 
\Benum
\Item  $\tau\subseteq \omega^2$ is \emph{a tree} iff
i)  for every $x\in Rng(\tau)$ the set $Pred(x):=\left\{y\in \omega;
\tau(y,x)\right\}$ is finite and $\tau\lceil Pred(x)$ is a linear
ordering and ii)
 there is $0_\tau\in \omega$ such that for every $x\in Rng(\tau)$,
$\tau(0_\tau,x)$.
\Item $\tau$ will also be written as $\preceq_\tau$.
$x\prec_\tau y$ is defined in the obvious way, and so is $inf(X)$ for
$X\subseteq Rng(\tau)$.

\noindent
$x\prec_\tau^s y$ iff $x\prec y$ and there is no $z$ $x\prec_\tau z$,
$z\prec_\tau y$.

\noindent
$Succ_\tau(x):=\{y;x\Preceq y\}$
$x\in T_\tau$. 

\noindent
$\tau_x$ is a tree such that $\tau_x:=\tau\lceil
Succ_\tau(x)$
\Item $B\subseteq T_\tau $ is \emph{a chain in $\tau$} iff
for every $x,y\in B$ $x\preceq y$ or $y\preceq x$.
\emph{ A branch in $\tau$ } is a maximum chain in $\tau$.
\Eenum
\Enddefinition

\Definition{defIIae} Let $\tau$ be a tree.
On $Rng(\tau)$ we define a binary relation $\K$ in the following way:

$x,y\in\omega$ then $\K(x,y)$ iff $y\Preceq x$ or there are $u,v\in Rng(\tau)$  $inf(x,y)\Precs
u,v$, $u\Preceq x$, $v\Preceq y$ and $u<v$. We shall refer to $\K$ as
\emph{ the Kleene-Brouwer ordering.}
\Enddefinition

The following two Propositions give us the basic properties of $\K$ that
we shall need. 
 
\Proposition{ThIIg} Let $\tau$ be a tree. Then
\Benum
\Item $\K$ is a linear ordering on $Rng(\tau)$.
\Item for every $x\in Rng(\tau)$, $KB(\tau_x)=\K\lceil
      Succ_\tau(x)$ and $Succ_\tau(x)$ is an interval in $\K$. 
\Item $\tau$ has no infinite branch iff $\K$ is a well-ordering.
\Item if $A=a_0,a_1,\dots $ is an infinite decreasing sequence in $\K$ then
the set $\{b_i; i\in\omega\}$, where $b_n:=
hinf\{ a_n,a_{n+1},a_{n+2}\dots \}$, is an infinite
chain in $\tau$.
\Eenum
\EndProposition

\Definition{defIIn}  Let $R$ be a linear ordering.
$wo(R)\subset Dom(R)$ is the set such that
 ${ wo}(R)$ is a maximum lower segment in $R$ such that $\prec_R\lceil
 wo(R)$ is a well-ordering.
$wo(R)$ will be called \emph{ the well-ordered part of $R$}.
The order-type of $\prec_R\lceil { wo}(R)$ shall be denoted by
      $ord(R)$.
\Enddefinition

It is evident that $wo(R)$ is defined uniquely; the existence follows
from the axiom of choice.

\Proposition{ThIIh}  \m, let $\tau$ be a tree definable in $\M$.
Then \Benum
\Item  $\K$ is definable $\M$.
In addition, if $\tau$ is $\Delta_1$ in $\M$ and for every $x,y$,
$x\Prec y$ implies $x<y$ then $\K$ is $\Delta_1$-definable in $\M$.
\Item If there is a nonempty $X\subseteq Dom(\K)$definable in $\M$ which does not have
      a $\K$-first member then there is an infinite branch of $\tau$
      definable in $\M$.
\Item Assume that $\tau$ has an infinite branch.
      Then if $wo(\K)$ is definable in $\M$ then an infinite branch of
      $\tau$ is
      definable in $\M$.
\Eenum
\EndProposition
\Proof 1) is obtained be translating the definition of $\K$ to the
structure $\M$.
2) and 3) follow from Proposition~\ref{ThIIg}. 
\Endproof

\bigskip
We shall now proceed to assign trees to formula schemes.
In the following definition, ${\mathbf y}$ will designate a list of
variables $y_1,\dots y_n$. $\forall \mathbf y$ will is an abbreviation
for   $\forall y_1,\dots \forall y_n$, and similarly in the case of
$\exists\mathbf y$.

\Definition{defIIaf}
Let $\psi$ be a closed formula scheme in a prenex form.
\Benum
\Item  
For $i=1,\dots n$ let 
 ${\mathbf y_i}={y^i}_1,\dots {y^i}_{l_k}$. We can assume that the
variables are mutually different.
Let 
\[
\psi=(\forall {\mathbf x_1}) (\exists {\mathbf y_1})\dots (\forall{\mathbf
  x_n})
(\exists {\mathbf y_n}) \psi_0
,\] 
where $\psi_0$ is an open formula. Let $k:=\sum_{i=1}^{n}l_i$.
Then \emph{ the Skolem formula for
$\psi$}, $\psi_s$ \rm will be the formula

$(H^1({\mathbf x_1},{y^1}_1) \land\dots  H^2({\mathbf x_1},{y^1}_{l_1})
 \land\dots
H^{k-l_n}({\mathbf x_n},{y^n}_1) \land\dots  H^k({\mathbf x_n},{y^n}_{l_n})
)\Rightarrow \psi_0$,

\noindent
 where $H^1,\dots H^k$ are second order variables not occuring in
$\psi$ of the appropriate arity.

\Item Let $\psi=\psi[X_1,\dots X_k]\in Fle^\lambda$, $\lambda\in \omega^{<\omega}$.
If $\psi_s$ is as above, we shall write $\psi_s=\psi_s[H^1,\dots
H^k]_s=\psi[H]=\psi_s[X_1,\dots X_k][H^1,\dots
H^k]_s$, where $H=\langle H^1,\dots
H^{k}\rangle$.     If $\mu=n_1,\dots n_k$ and the arity of $H^i$ is
$n_i$ then we shall say that
$\psi_s\in Fle^{\lambda,\mu}$.

\Item Let $\psi$ be a formula in a prenex form,
$\psi_s=\psi_s[H^1,\dots H^k]_s$, where the arity of $H^i$ is $n_i+1$.
Functions $g_i:\omega^{n_i}:\rightarrow\omega$, $i=1,\dots k$, will
be called the 
\emph{ Skolem functions for
$\psi$  } iff $\psi_s[~^{n_1+1}g_1,\dots ~^{n_k+1}g_k]$ is true.
\Eenum
\Enddefinition

In the third item of the definition we identify $n$-ary function
$g:\omega^n\rightarrow \omega$ with the set of $(n+1)$-tuples and hence we can use a
function in the place of a predicate. Note that $n$-ary function $f$ can stand in a place of
a function symbol as $^nf$, while in a place of a predicate as
$^{n+1}f$.
The following is then obvious:
\Smallskip \it 
 Let
 $\psi=\psi(y_1,\dots y_{n},z)$ be a formula scheme. Let $f:\omega^n\rightarrow
 \omega$.
Then the scheme  $\psi(y_1,\dots y_n,~^nf(y_1,\dots y_n))$ is equivalent
to the scheme
$\forall z~^{n+1}f(y_1,\dots y_n,z)\Rightarrow\psi(y_1,\dots
y_{n},z)$.
\Smallskip\rm
We may conclude that if $\psi_s=(H^1({\mathbf y_1},z_1)\land
\dots H^k({\mathbf y_k},z_k))\Rightarrow \psi_0$ then
the functions
$g_1,\dots g_k$ 
 are Skolem functions for $\psi$ iff the formula
$\forall {\mathbf y_1}\dots \forall{\mathbf y_k}\psi_0(z_1/~^{n_1}g_1({\mathbf y_1}),\dots
z_k/~^{n_k}g_k({\mathbf y_k}) )$ is true. This implies the following lemma:

\Lemma{LemmaIIan} \m. Let $\lambda,~\mu\in \omega^{<\omega}$, let $\psi\in
Fle^\lambda$ be a formula scheme in
$\M$ in a prenex form. Let $B\in P(\omega^\lambda)$. Let $\psi_s\in
Fle^{\lambda,\mu}$. Then
\Benum
\Item  $B$ satisfies $\psi$ iff there are functions $g_1,\dots g_l$
which are the Skolem functions for $\psi_s[~^\lambda B]$.

\Item Assume that $B$ is definable in $\M$ and $g_1,\dots g_l $ are
Skolem functions for $\psi_s[~^\lambda B]$. Let  $h_i:=g_i\lceil  n,~
i=1,\dots l$. Then there are Skolem functions $w_1,\dots w_l$ for $\psi_s[B]$
such that  $h_i\subseteq
w_i,i=1,\dots l$ and $w_1,\dots
w_l$ are definable in $\M$.
\Eenum
\EndLemma

\Definition{defIIq} Let $\lambda,\mu\in \omega^{<\omega}$, $\lambda=\langle
a_1,\dots a_s\rangle$, $\mu=\langle b_1+1,\dots b_t+1\rangle$.
 Let $\psi\in Fle^\lambda$ be a formula scheme
in a prenex form.
Let $\psi_s=\psi_s(y_1,\dots y_l)\in Fle^{\lambda,\mu}$.
Then $\alpha$ is a \emph{ satisfaction system of degree
$n$ for $\psi$ }  iff
\[
\alpha =\langle B_1\dots B_s,g_1\dots g_t \rangle
\] 
and the following
conditions are satisfied\footnote{Recall that $n=\left\{0,1\dots
n-1\right\}$, $n\in \omega$.}.
\Benum
\Item $g_i: n^{b_i}\rightarrow \omega, i=1\dots t$

\Item Let $e_\alpha:=\max\left\{n,\max(Rng (g_i))+1;i=1\dots t\right\}$.
Let $F$ be the set of functions occuring in $\psi$. Let
$m_\alpha:=\max\left\{n,
\max Rng(f\lceil e_\alpha);i=1\dots t, f\in F\right\}$.
Then $B_i\subseteq m_\alpha^{a_i}$, $i=1\dots s$

\Item The formula $\psi_s[~^\lambda B][~^\mu g](y_1,\dots y_l)$ is
      satisfied by every $a_1,\dots a_l$ such that $a_1,\dots a_l<e_\alpha$.
\Eenum
\Enddefinition

The intuition behind the definition is simple.
Assume, for clarity, that $\psi_s$ contains
no function symbols and that $\psi_s=\psi=\psi[A_1,\dots A_n]$
(i.e. $\psi$ is an open formula scheme containing no
function symbols). Then a satisfaction system of degree $n$ is a sequence of
$B_i\subseteq\{ 0,1,\dots n-1 \}^{a_i}$, $i=1,\dots s$, such that the formula
$\psi_s[~^{b_1}B_1,\dots~^{b_s} B_s]$ is true when we let the variables range over
$0,1,\dots n-1$ only. The sets $B_i$ can be viewed as predicates 
defined on $\{ 0,1\dots n-1\}^{b_i}$, and we demand they satisfy the
formula on the domain of their definition. In general $\psi_s=\psi_s[A_1,\dots A_n][g_1,\dots g_l]$ and we assume that
the functions $g_1,\dots g_l$ are defined on $\{ 0,1,\dots n-1\}$. However,
we must make sure that the predicates $B_1,\dots B_s$ are defined on the ranges
of those functions and the other functions occuring in $\psi_s$ restricted
on $\{ 0,1,\dots n-1\}$. Point 2) reflects the fact that for a Skolem function $f$
and $g\in{\cal F}(L)$ the term $g(f)$ may occur in $\psi_s$ but the term $f(g)$ cannot.

\Definition{defIIr} Let $\psi\in Fle^\lambda$ be a formula scheme  in a prenex
form. 
\Benum
\Item
Let $\alpha=\langle B_1\dots B_s,g_1\dots g_t \rangle$ be a
satisfaction system for $\psi$ of degree $n$ and let
$\alpha^\prime=\langle B_1^\prime\dots 
B_s^\prime,g_1^\prime\dots g_t^\prime \rangle$ be a
satisfaction system for $\psi$ of degree $m$. Then we let 
\[
\alpha\preceq\beta
\]
iff $m\leq n$ and i) for every $i=1,\dots t$ if $g_i$ is a $k$-ary
function then $g_i=g_i^\prime\lceil n^k$
ii) for every $i=1,\dots s$ if $B_i\subseteq \omega^k$ 
then $B_i=B_i^\prime\cap e_n^k$.
\Item \emph{ The characteristic tree of $\psi$ } is the $\tau\subseteq
\omega^2$ such that
for every $a,b\in \omega$, $\tau(a,b)$ iff 
there are $\alpha,\beta$ satisfaction systems of $\psi$,
$a=[\alpha],~b=[\beta]$\footnote{Here $[\alpha]$ is the G\"odel number of
the finite set $\alpha$ see page \pageref{code}.}
 and $\alpha\preceq\beta$.
\Item Let  $\psi$ be a formula scheme in a
prenex form. Let $\tau_\psi$ be the characteristic tree of $\psi$.
 The ordinal $ord(KB(\tau_\psi))$ will be called
\emph{ the characteristic ordinal of $\psi$ }; it will be denoted by
$ord(\psi$).
\Eenum
\Enddefinition

The key property of a satisfaction system is expressed in the next
Lemma.

\Lemma{LemmaIIo} Let  $\psi\in Fle^\lambda$ be in a prenex form. Let
 $\alpha^i=\langle B_1^i\dots 
B_s^i,g_1^i\dots g_t^i \rangle$, $i\in \omega$ be a sequence of
satisfaction systems of $\psi$ such that
$\alpha_0\prec\alpha_1\prec\alpha_2\dots$. Then
$B:=\langle\bigcup_{i\in\omega} B^i_1,\dots 
\bigcup_{i\in\omega} B^i_s\rangle$ satisfies $\psi$ and
$g_l:=\bigcup_{i\in\omega}g_l^i,~l=1,\dots t$, are the Skolem functions
for $\psi[~^\lambda B]$.
\EndLemma
\Proof Evident
\Endproof

\Proposition{ThIIaf}  Let $\M$ be a structure. Let $\psi$ be a formula scheme
in a prenex form in $\M$, $\tau_\psi$ the characteristic tree of
$\psi$. Let $\cal A$ be the system defined by $\psi$.
 Then
 \Benum
 \Item  $\tau$ is a $\Delta_1$-definable tree in $\M$.
 
\Item  
$\cal A\not=\emptyset$ iff $\tau$ contains an infinite branch.

\Item There is $B\in \cal A$ definable in $\M$ iff $\tau$ contains a definable
infinite branch in $\M$.
\Item $KB(\tau_\psi)$ is $\Delta_1$-definable in $\M$.
      If $\cal A\not=\emptyset$ and $wo(KB(\tau_\psi))$ is definable in
      $\M$ then some $A\in \cal A$ is definable in $\M$.
\Eenum
\EndProposition

\Proof 
1) follows from   the definition of $\tau_\psi$ and  Proposition~\ref{Thc}. 

2) and 3) follow from  Lemmata \ref{LemmaIIo} and \ref{LemmaIIan}.

4) follows from Proposition~\ref{ThIIh}, 1) and 3).
\Endproof

\Lemma{LemmaIIp} \m, $\psi$ a formula scheme in $\M$ in a prenex form. Let $\psi$ be a
formula scheme in $\M$ which defines a non-empty system $\cal A$. Let $\ro$ be an arbitrary representation of ordinal
$\alpha$, $ord({\cal A})\leq\alpha$.  Let $c$ be a coding for the
structure $\M$ 
compatible with $\ro$. Then there is $B\in {\cal A}$ definable in
$\A+~^2\Trro$.
\EndLemma
\Proof Let $\tau$ be the characteristic tree of $\psi$.
By Theorem~\ref{ThIIaf}, 4) it is sufficient to show that ${ wo}(\K)$ is
definable in $\A+~^2\Trro$

Since $\K$ is not a well-ordering, then $wo(\K)\not=\omega$. Let us
chose $a\not\in wo(\K)$.

Let $f$ be a function $f:Dom(\ro)\rightarrow Dom(\K)$ such that
\Benum 
\Item if there is $z$ a $\K$-minimum of $Rng(\tau)$, let $f(\ro_0):=z$ else
$f(\ro_0):=a$ 
\Item Let $y\in Rng(\tau)\setminus \{\ro_0\}$.  If there exists $z$
which is $\K$-minimum of
$Rng(\K)\setminus Rng(f\lceil \left\{x;x\prec_\ro y\right\} ) $,
then $f(y):=z$.
Otherwise $f(y)=a$.
\Eenum
Since $\K$ is definable in $\M$, it is trivial to find $B$ and $F$
definable in $\M$ such that 
\[
f=REK_2(B,F,\ro)
\]
We can apply Proposition~\ref{ThIIe} to obtain that $f$ is definable in $\A+~^2\Trro$.
But $wo(\K)=Rng(f)\setminus \left\{a\right\}$ and hence $wo(\K)$ is
definable in $\A+~^2\Trro$.
\Endproof

\Lemma{LemmaIIq} \m.
Let $R$ be a linear ordering definable in $R$.
Then \[
ord(R)\leq Ord(\M)
\]
 If in addition $R$ is $\Delta_1$ in $\M$ then $ord(R)\leq\xi$, the
 first ordinal not $\Delta_1$-definable in $\M$. 
\EndLemma
\Proof Let $R$ be a linear ordering definable in $\M$.  For $n\in {wo}(R)$,    $R_n$
is a representation of an ordinal $\alpha_n\leq ord(R)$\footnote{For
  the definition  of $R_n$ see page \pageref{Rn}}.
 Clearly $ord(R)=sup \left\{ \alpha_i,i\in \omega\right\}$.
If $R$ is ($\Delta_1$-) definable in $\M$ then $R_n$ is ($\Delta_1$-)
definable in $\M$ for every $n\in wo(R)$. Hence $ord(R)\leq Ord(\M)$ (resp. $ord(R)\leq \xi$).
\Endproof

 \Proposition{ThII}
\m. Let $\alpha$ be a countable ordinal.
Then the following conditions are equivalent
\Benum
\Item$\alpha$ is implicitly definable in $\M$.
       
\Item$\alpha$ is definable in $\M$.
                
\Item$\alpha$ is $\Delta_1$-definable in $\M$.
\Eenum
\EndProposition
\Proof The implications 3)$\Rightarrow$ 2) $\Rightarrow$ 1) are trivial.

1)$\Rightarrow$ 3). Assume the contrary.
Without the loss of generality we can assume that $\M$ is finite.
 Let $\ro$ be an implicitly definable
representation of $\alpha$ such that $\alpha$ is not
$\Delta_1$-definable. 
We can assume that $\ro$ is compatible with a coding $c$ for $\M$. We
will show that ${\cal I}(\M)\subseteq \D(\A+^2\Trro)$.
Let $\psi$ be a formula scheme in $\M$ which
defines a non-empty system $\cal A$. By the Lemma~\ref{LemmaIIq} every
$\beta<ord(\psi)$ is $\Delta_1$-definable in $\M$. It follows that
$\alpha\geq ord(\psi)$. By Lemma~\ref{LemmaIIp} there is some $A\in \cal A$ definable in
$\A+~^2\Trro$.
Hence ${\cal I}(\M)\subseteq \D(\A+^2\Trro)$. We assumed that $\M$ is finite
and  $\ro$ is implictly definable in $\M$  and hence also $\Trro$
is implicitly definable in $\M$; therefore  ${\cal I}(\M)\sim
\A+^2\Trro$. But that contradicts Corollary 2 of Proposition~\ref{ThIIa}.
\Endproof

\bigskip
The
statement of the proposition may be strengthened to say that every
definable system of ordinals contains a definable element or even that
every definable system of ordinals has a definable supremum (see [2],
Chapter IV), but those modifications 
will not be needed here.

\Proposition{ThIIad}\m. Let 
$\psi$ a formula scheme in $\M$. Then
\Benum
\Item   $ord({\psi})\leq Ord(\M)$ 

\Item If $\psi$ is a proper implicit definition of some $B$ then
$ord(\psi)< Ord(\M)$.
\Eenum
\EndProposition

\Proof 1) Let $\tau$ be the characteristic tree of $\psi$ in $L$ which
defines $\psi$. Then $\K$ is
definable in $\M$ (Proposition~\ref{ThIIh}) and $ord(\K)\leq Ord(\M)$ by
Lemma~\ref{LemmaIIq}.

2) without the loss of generality we can assume that $\M$ is finite
and that $B\subseteq \omega$.
It is sufficient to prove that $wo(\K)$ is implicitly definable in
$\M$. Since $\K$ is definable in $\M$, $\K\lceil wo(\K)$ is then an
implicitly definable
representation of $ord(\psi)$ in $\M$ and therefore, by Proposition~\ref{ThII},
$ord(\psi)$ is definable in $\M$. Let $T$ be a truth predicate for
$\M+~^1B$. Since $T$ is implicitly definable in $\M+~^1B$ and $B$ is
implicitly definable in $\M$, $T$ is implicitly definable in $\M$ and
it is sufficient to prove that $wo(\K)$ is definable in
$\M+~^1B+~^1T$.

Clearly, the two conditions are equivalent: \Benum
\Item $x\in wo(\K)$
\Item there is no $y\in Rng(\tau)$ such that $\K(y,x)$ and $y$ lies on an infinite
      branch of $\tau$.
\Eenum  
On the other hand, by Lemma~\ref{LemmaIIan}, 2) the condition `$y$ lies on an
infinite branch of $\tau$` is equivalent
to the condition `$y$ lies on an infinite branch of $\tau$ definable
in $\M+~^1B$`. But the later statement can be expressed using the truth
predicate $T$, and hence the condition 2) can be expressed in $\M+~^1B+~^1T$.
Therefore $wo(\K)$ is definable in $\M+~^1B+~^1T$.
\Endproof

\Lemma{ThIIaa} Let $\N$ be a structure.
Let ${\cal T}$ be a set of binary predicates such that
i) every $X\in \cal T$ is of the form $~^2Tr_{\N,c,\ro}$, where $\ro$ is
a well-ordering definable in $\N$ compatible with coding $c$ and 
ii) for every $\alpha<Ord(\M)$ there exists $\ro$
a representation of $\beta\geq\alpha$ definable in $\N$ and a coding $c$ for
$\N$ such that 
$~^2Tr_{\N,c,\ro}\in \cal T$. Then 
${\cal I}(\J)\sim A+{\cal T}$.
 
\EndLemma
\Proof ${\cal I}(\J)\subseteq \D(A+{\cal T})$  follows from Lemma~\ref{LemmaIIp} and Proposition~\ref{ThIIad}.
By Proposition~\ref{ThIIae} we have also ${\cal T\subseteq \cal
  I}(\N)$
\Endproof

\bigskip
\noindent
{\tt \bf Proof of Theorem~\ref{Thx}.}\label{proofx}
Theorem~\ref{Thx} is now a direct consequence of the previous Lemma and the corollary of Theorem~\ref{LemmaIIIoq}.
\Endproof
\bigskip

In order to prove Theorem~\ref{Thy}, we shall find a linear ordering $R$ definable in
$\N$ such that $ord(R)=Ord(\N)$. This will be achieved by means of a
formula scheme in $\N$ such that $ord(\psi)=Ord(\N)$.
It must be noted that  $ord(\psi)\leq Ord(\J)$ but
 the condition $ord(\psi)<Ord(\J)$
in general holds just for schemes which implicitly define a set. In $\N$
there may exist systems defined by a scheme $\psi$ such that the
characteristic ordinal of $\psi$ is not definable in $\N$. 
 
Observe  that for a formula $\psi$ defining a
nonempty system $\cal A$ if $ord({{\psi}})<Ord(\M)$ then, by
Lemma~\ref{LemmaIIp}, there is some $A\in {\cal A}$
implicitly definable in $\M$. Hence, if for
every formula scheme in $\M$, $ord({{\psi}})<Ord(\M)$ then 
every non-empty system
definable contains an implicitly definable set in $\M$. This is the
essence of the folowing definition.

\Definition{defIIah } Let $\M$ be a structure.
Then  $\M$ is \emph{ implicitly complete} iff
every non-empty system definable in $\M$ contains an implicitly definable set in
$\M$. 
\Enddefinition

\Lemma{ThIIn}  Let $\M$ be a structure.
If for every $R$ a linear ordering definable in $\M$,
$ord(R)<Ord(\M)$ then $\M$ is implicitly complete.
\EndLemma

\Proof We can assume that $\M$ is finite.
Let $\psi$ be a
scheme in $\M$ in a prenex form which
defines a non-empty system ${{\cal A}}$. Let $\tau$ be the
characteristic tree. $\K$ is definable in $\M$ and by the assumption there
is $ord(\K)<Ord(\M)$. We can find a definable representation $\ro$ in
$\M$ of an ordinal $\beta$, $ord(\K)<\beta<Ord(\M)$. We can assume
that $\ro$ is compatible with a coding $c$ for $\M$. The set $\Trro$
is implicitly definable in $\M$ and by Lemma~\ref{LemmaIIp} there is
some $A\in \cal A$ definable in $\A+~^2\Trro$. 
\Endproof

\bigskip
In definition \ref{Def4a2} we introduced $\overline{Tr}(\M,c,R)$,which
is a
generalisation of the concept of $\Trro$ if $R$ is not a
well-ordering. Similarly to Proposition~\ref{ThIIae} we may obtain: 
\Smallskip \it
There is a system ${\cal S}\subseteq P(\omega^{2,2,2})$ definable in
$\A$ such that for
every structure $\M$ and a coding $c$ for $\M$ and a linear ordering
$R$,  $\langle D_{\M,c}, R, U\rangle \in \cal S$ iff 
$U\in\overline{Tr}(\M,c,R)$
\Smallskip \rm

\Lemma{LemmaIIag}  Let $\J$ be a $L$-finite structure, $c$ a coding for $\J$.
Let $R$ be a linear ordering compatible with $c$ such that $Ord(\J)\leq ord(R)$.
Let $U\in \overline{Tr}(\J,c,R)$.
If  $x\in Rng(R)\setminus wo(R)$ then there is a
universal predicate for ${\cal I}(\J)$ definable in
$A+~^1U(x,.)$. Hence ${\cal I}(\N)\subseteq \D(\A+~^1U(x,.))\subseteq  \D(\A+~^2U)$.
\EndLemma
\Proof  For $n\in N$ let $U_n\subseteq N^2$ be a relation such that  $U_n(a,b)$
iff $U(a,b)$ and $R(a,n)$, let $R_n$ be defined as on page  \pageref{Rn}. For  $n\in wo(R)$ we have $U_n=Tr_{\N,c,R_n}$.
Since $Ord(\N)\leq ord(R)$ then for every
ordinal $\alpha$ definable in $\J$ there is some $n\in wo(R)$ such that $R_n$ is
a representation for $\alpha$. 
Therefore, using Lemma~\ref{ThIIaa}, every set implicitly definable in $\J$ is
definable in $\J+\{~^2U_n;n\in wo(\N)\}\sim \J+\{~^1U(n,.);n\in wo(\N)\}$ (see Lemma~\ref{LemmaIIac}).
The set $U(x,.)$ is defined to be a truth predicate for the structure
$\M:=\A+\{ ~^1U(y,.);{y\prec_Rx)}\}$. A proper universal predicate $P$ for $\M$ is
therefore definable in $\A+~^1U(x,.)$. But 
${\cal I}(J)\subseteq \D( \J+\{~^2U(n,.);n\in
  wo(\N)\})\subseteq \M$ and hence $P$ is a universal
predicate for ${\cal I}(\J)$.
\Endproof

\Proposition{ThIIp} Let $\J$ be a $L$-finite structure.
Then $\J$ is not implicitly complete.
\EndProposition

\Proof Let $\J$ be a $L$-finite structure and $c$ a coding for $\J$. In $\J$ we
can find a formula scheme in $\psi\in Fle^{2,2,1,3}_1(\J)$ such that for
every $R,U,X,Y\in P(\omega^{2,2,1,3})$, $R,U,X,Y$ satisfies $\psi$ iff
   
\Benum
\Item $R$ is a linear ordering compatible with $c$ and $R$ has the
      first  member $0_R$,
\Item $U\in \overline{Tr}(D_{\N,c},R)$,

\Item for every $n\in N$ if $n$ is a $c$-G\"odel
      number of a formula in $\J$ defining a linear
      ordering $Q$ then either a)
$X(n,.)\subseteq Dom(Q)$ is a set without the $Q$-first member or b)
$Y(n,.,.)$ is an isomorphism between $Q$ and a lower segment of $R$.
\footnote{I.e., for every $x\in Dom(Q)$ there is a unique $y\in Dom(R)$ such that
$X(n,x,y)$ and i) if $0_Q$ is the $Q$-first member of $Dom(Q)$ then
$Y(n,0_Q,0_R)$ and ii) for every $x\in Dom(Q)$ $Y(n,x,y)$ iff $y$ is
the $R$-first member of $Rng(R)\setminus \{z;\mbox{ there exists }z^\prime\in
Dom(Q),~ 
Y(n,z^\prime,z)~and~z^\prime\prec_Q x\}$.}  
\Eenum
Observe that the condition 3) uses just quantifications over sets
definable in $\J$ and hence it can be formulated using the fact that $R$ has the
smallest member and $U(R_0,.)=\Tr$.

If $R,U,X,Y$ satisfy $\psi$ then $Ord(\J)\leq ord(R)$ because for every
well-ordering $Q$ definable in $\J$ the condition a) is not satisfied and
therefore $Q$ must be isomorphic to a lower segment of $R$.
Moreover, for
every  $R$ which is a representation of an ordinal $\alpha\geq Ord(\J)$
compatible with $c$ there are some $U,X,Y$ such that $R,X,Y,U$ satisfies
$\psi$. Hence the system defined by $\psi$ is non-empty.

Let us assume that $\J$ is implicitly complete. Then there are some
$R,X,Y,U$ satisfying $\psi$ which are implicitly definable in $\J$.
By Lemma~\ref{LemmaIIag} we have  ${\cal I}(\J)\subseteq \D(\A+~^2U) $.
Since we assumed that $U$ is implicitly
definable this implies that the structure ${\cal I}(\J)$ is
essentially finite. But this contradicts the Corollary 2 of Proposition~\ref{ThIIa}.
\Endproof

\bigskip
\Corollary{} Let $\N$ be finite. Then there is a linear ordering $R$ definable in
$\N$ such that
$ord(R)=Ord(\J)$. Hence $ord(R)$ is not definable and $wo(R)$ is not
implicitly definable in $\M$.
\Endcorollary

\Proof Follows from the previous Theorem and Lemma~\ref{ThIIn}.
\Endproof.

\Lemma{Lemma4a2} Let $\N$ be a $L$-finite structure and $c$ a coding for
$\N$.
 Then there is a linear
ordering $R$ definable in $\N$ and some $U\subseteq \omega^2$ with the following properties:
\Benum
\Item $ord(R)=Ord(\M)$,
\Item $R$ is compatible with $c$ and $U\in \overline{Tr}(\M,c,R)$,
\Item there is no $X\subseteq Rng(R)$ such that $X$ does not have  the $R$-first member
      and $X$ is definable in $\N+~^2U$.
\Eenum
\EndLemma

\Proof Assume the contrary. Let $R$ be a definable linear ordering in
$\N$ such that $ord(R)=Ord(\N)$. Let $c$ be a coding for $\N$.
We can assume that $R$ is compatible with $c$.
 Let $\psi\in Fle^{1,1,2}(\N)$ be a formula scheme such that $T,V,U$ satisfy $\psi$ iff
\Benum
\Item $V\subseteq Rng(R)$ is a lower segment in $R$ such that $Rng(R)\setminus V$ is
      non-empty and does not have a $R$-first member.
\Item $U\in \overline{Tr}(\N,c,R\lceil V)$.
\Item Let $c^\prime$ be the coding for the structure $\N+~^2U$ such
      that $c\subseteq c^\prime$ and
      $c^\prime(~^2U)=min(\omega\setminus Rng(c))$.
      Then $T=Tr_{\N+~^2U,c^\prime}$. 
\Item There is no non-empty $X\subseteq Rng(R)$ definable in
      $\N+~^2U$ such that $X$ does not have a $R$-first element.
\Eenum
Observe that the last condition can be formulated using the truth
predicate $T$.

Let us show that under the given assumption the formula $\psi$
is a proper implicit definition of some $T,V,U$ such that
$V=wo(R)$.

Assume first that $V=wo(R)$. 
Then there is unique $U$ which satisfies 2) because
$R\lceil V$ is a well-ordering. Then there is unique $T$ such that
3)  is satisfied. The condition 4) is satisfied because  
$R\lceil V$ is a well-ordering.

Assume that $V\subseteq wo(R)$, $V\ne wo(R)$. Then clearly 1) is not
satisfied.

Assume that $wo(R)\subseteq V$, $V\ne wo(R)$. Then 2) or 4)  is not satisfied by
the assumption.

Hence  $\psi$ is an implicit definition of some $T,wo(R),U$.
But that is impossible. For then $R\lceil wo(R)$
is implicitly definable representation of $Ord(\N)$ and
hence $Ord(\N)$ is implicitly definable, contrary to Proposition~\ref{ThII}.
\Endproof

\bigskip
\noindent
{\tt \bf Proof of Theorem~\ref{Thy}\label{proofy}.} 
Assume that $\cal I(N)$ has a canonic universal predicate $G_0$.
Let us first prove the following:  
\Smallskip
\it Let $R$ be  a linear ordering
definable in $\N$ and let $U\in
\overline{Tr}(\N,c,R)$ for a compatible coding $c$ for $\N$. Assume that $U$ is definable in $\N+~^2G_0$. Then 
$ord(R)<Ord(\N)$ \rm
\Smallskip
Let $R,U$ be as assumed and $ord(R)=Ord(\N)$. 
Let $m,n\in
Rng(R)\setminus wo(R)$, $m\prec_R n$. (The existence is granted since
$wo(R)$ is not definable in $\N$ and hence $Rng(R)\setminus wo(R)$
must be infinite). Let $U_m:=U(m,.)$ and $U_n:=U(n,.)$. By Lemma
\ref{LemmaIIag} a universal predicate for $\cal I(N)$ is definable in
$\N+~^2U_m$. 
Therefore also $G_0$ is definable in $\N+~^2U_m$ because $G_0$ is
canonic. 
But since $U$ is definable in $\N+~^2G_0$ then $U_n$ is definable in
$\N+~^2G_0$ and hence also in $\N+~^2U_m$. 
But from the definition of $U\in \overline{Tr}(\M,c,R)$, $U_n$ is a
truth predicate for a structure containing $\N+~^2U_m$. But that is
impossible.

 Let us now complete the proof. Let us take $R,U$ as in the previous
 Lemma. By Lemma~\ref{LemmaIIag} there is  a universal predicate for $\cal
 I(M)$ definable in $\N+~^2U$. By the assumption, $G_0$ is definable in
 $\N+~^2U$. 

In  $\N+~^2G_0$ we can find a formula $\eta$ with one free
 variable such that:
\Smallskip\it
 for every $n\in \omega$
$n$ satisfies $\eta$ iff there exists $m\in \omega$ such that $m$ is a
$G_0$-code of some $U_n$ such that $U_n\in \overline{Tr}(\N,c,R_n)$ 
\Smallskip\rm
If $n\in wo(R)$ then $R_n$ is a well-ordering definable in $\N$. Hence
$Tr_{\N,c,R_n}$ is implicitly definable in $\N$ and it has a
$G_0$-code since $G_0$ is a universal predicate for $\cal I(N)$. Hence $n$
satisfies $\eta$. If on the other hand $n\not\in wo(R)$ then $n$ does
not satisfy $\eta$ by the proposition. This implies that $wo(R)$ is
definable in  $\N+~^2G_0$ and therefore also in  $\N+~^2U$. But this
contradicts the condition 4) of the Lemma~\ref{Lemma4a2}, since $Rng(R)\setminus
wo(R)$ does not have a $R$-first member.
\Endproof

\bigskip
Recall the relation between  Tarski and proper Tarski hierarchy as
stated in Corollary~2 of Theorem~\ref{LemmaIIIoq}.
Hence, in order to prove  Theorem~\ref{Thb}, it is
sufficient to show that the structure $\M_{Ord(\N)}\sim
\M_{Ord(\N)}^p\sim {\cal I(N)}$ does have a canonic proper universal
predicate.

\Theorem{ThIIIt} Let $\J$ be a $L$-finite structure. Let
$\{ \M^p_\alpha\}_{\alpha<\lambda^p(\J)}$ be a proper Tarski hierarchy over $\J$.
Then
$Ord(\J)<\lambda^p(\J)$.
\EndTheorem

\Proof Let us show that the
structure $\M_{_{Ord(\J)}}$ does have a canonic proper universal predicate. 
 By Proposition~\ref{ThIIp} and Lemma~\ref{ThIIn} there is a linear ordering $R$
definable in $\J$ such that $ord(R)=Ord(\J)$.  Let
$\ro=R\lceil wo(R)$. 
$\omega\setminus wo(R)$ is
infinite and we can chose a coding $c$ for $\J$ compatible with $\ro$.
Furthermore,
for every $\gamma<Ord(\J)$ $\ro\lceil \gamma$ is a definable representation of
$\gamma$ in $\J$. Hence, from Theorem~\ref{LemmaIIIoq}, 
$$
\M_{_{Ord(\J)}}\sim\J+\{Tr_{\J,c,\ro}(\ro_\gamma,.);{\gamma<Ord(\J)}\}
$$
Let $\M_{<Ord(\J)}$ denote the structure on
the right hand side of the equivalence, and $c_{_{Ord(\J)}}$ be the induced
coding on $\M_{<Ord(\J)}$. As in the proof of Theorem~\ref{LemmaIIIoq} it is sufficient
to prove that  $Tr_{\M_{<Ord(\J)}, c_{_{Ord(\J)}}}$ is definable in
$\A+P$, for any proper universal predicate $P$ for
$\M_{<Ord(\J)}$. 

Let $P$ be a proper universal predicate for $\M_{<Ord(J)}$.
From  Lemma~\ref{LemmaIIIos} it is sufficient to prove that
$\ro$ is definable in $\A+P$. 
In  $\N+P$ we can find a formula $\eta$ with one free
 variable such that:
\Smallskip\it
 for every $n\in \omega$
$n$ satisfies $\eta$ iff there exists $m\in \omega$ such that $m$ is a
$P$-code of some $U_n$ such that $U_n\in \overline{Tr}(\N,c,R_n)$ 
\Smallskip\rm
If $n\in wo(R)$ then $R_n$ is a well-ordering definable in $\N$. Hence
$Tr_{\N,c,R_n}$ is implicitly definable in $\N$ and it has a
$P$-code since $P$ is a universal predicate for $\cal I(N)$. Hence $n$
satisfies $\eta$. If on the other hand $n\not\in wo(R)$ then 
there is no $U_n$ implicitly definable in $\N$ such that $U_n\in
\overline{Tr}(\N,c,R_n)$ (for otherwise ${\cal I(M)}\sim \A+~^2U_n$
and $\cal I(M)$ is essentially finite). Hence $n$ does
not satisfy $\eta$ because $P$ is a proper universal predicate for
$\cal I(M)$.  Therefore $wo(R)$ is
definable in $\A+P$. 
\Endproof

\bigskip
\Corollary{} 
There is a structure which has a canonic proper universal
predicate but does not have a canonic universal predicate.
Namely, if $\N$ is a $L$-finite structure then $\cal I(N)$ has a canonic
universal but not a proper canonic universal predicate.
\Endcorollary


\end{document}